\documentclass{gtart_h}

\def\ifplaintex{\expandafter\ifx\csname documentclass\endcsname\relax}

\def\gtp{{\mathsurround=0pt\it $\cal G\mskip-2mu$eometry \&\ 
$\cal T\!\!$opology $\cal P\!$ublications}}  

\def\recd{{\small Received:\qua\receiveddate\ifx\reviseddate\relax
\else\qquad Revised:\qua\reviseddate\fi\par}} 

\def\lognumber#1{\def\thelognumber{#1}}
\def\volumenumber#1{\def\thevolumenumber{#1}}
\def\volumeyear#1{\def\thevolumeyear{#1}}
\def\papernumber#1{\def\thepapernumber{#1}}
\def\pagenumbers#1#2{\def\startpage{#1}\def\finishpage{#2}}
\def\published#1{\def\publishdate{#1}}

\def\received#1{\def\receiveddate{#1}}
\def\revised#1{\def\reviseddate{#1}}
\def\accepted#1{\def\accepteddate{#1}}

\def\asciiurl#1{\def\theasciiurl{#1}}

\long\def\asciiabstract#1{\long\def\theasciiabstract{#1}}
\def\asciikeywords#1{\def\theasciikeywords{#1}}

\let\\\par\let\thelognumber\relax\let\thevolumenumber\relax
\let\thepapernumber\relax\let\thevolumeyear\relax\let\startpage\relax
\let\finishpage\relax\let\publishdate\relax\let\receiveddate\relax
\let\reviseddate\relax\let\accepteddate\relax\let\theasciititle\relax
\let\theasciiauthors\relax
\let\theasciiabstract\relax\let\theasciikeywords\relax

\let\theasciiemail\relax
\let\theasciiurl\relax

\ifplaintex
\font\logobig=cmssbx10 scaled 3836
\font\logomed=cmssbx10 scaled 2557
\else
\font\logobig=cmssbx10 scaled 4200
\font\logomed=cmssbx10 scaled 2800
\fi

\long\def\makeagttitle{   
\count0=\startpage
\agt\hfill      
\hbox to 45truept{\vbox to 0pt{\vglue -13truept{\logomed A\kern -.37em{\logobig 
T}\kern -.38em G}\vss}\hss}
\break
{\small Volume \thevolumenumber\ (\thevolumeyear)
\startpage--\finishpage\nl
Published: \publishdate}

\vglue .25truein

{\parskip=0pt\leftskip 0pt plus
1fil\def\\{\par\smallskip}{\Large\bf\thetitle}\par\medskip} \vglue
0.05truein

{\parskip=0pt\leftskip 0pt plus 1fil\def\\{\par}{\sc\theauthors}
\par\medskip}%
 
\vglue 0.03truein 

{\small\leftskip 25truept\rightskip 25truept{\bf Abstract}\stdspace\theabstract

{\bf AMS Classification}\stdspace\theprimaryclass
\ifx\thesecondaryclass\relax\else; \thesecondaryclass\fi\par
{\bf Keywords}\stdspace \thekeywords\par}\vglue 7truept

}   

\ifplaintex
\font\phead=cmsl9 scaled 950
\font\pnum=cmbx10 scaled 913
\font\pfoot=cmsl9 scaled 950
\headline{\vbox to 0pt{\vskip -4.5mm\line{\small\phead\ifnum
\count0=\startpage ISSN 1472-2739 (on-line) 1472-2747 (printed)
\hfill {\pnum\folio}\else\ifodd\count0\def\\{ }%
\ifx\theshorttitle\relax\thetitle\else\theshorttitle\fi\hfill{\pnum\folio}
\else\def\\{ and }{\pnum\folio}\hfill\ifx\theshortauthors\relax\theauthors
\else\theshortauthors\fi\fi\fi}\vss}}
\footline{\vbox to 0pt{\vglue 0mm\line{\small\pfoot\ifnum\count0=\startpage
\copyright\ \gtp\hfill\else
\agt, Volume \thevolumenumber\ (\thevolumeyear)\hfill\fi}\vss}}
\else
\font=cmsl9 scaled 1050
\font\lnum=cmbx10 
\font=cmsl9 scaled 1050
\makeatletter
\def\@oddhead{{\small\lhead\ifnum\count0=\startpage ISSN 1472-2739 
(on-line) 1472-2747 (printed)\hfill {\lnum\number\count0}\else\ifodd\count0
\def\\{ }\ifx\theshorttitle\relax \thetitle \else\theshorttitle\fi\hfill
{\lnum\number\count0}\else\def\\{ and }{\lnum\number\count0}
\hfill\ifx\theshortauthors\relax 
\theauthors\else\theshortauthors\fi\fi\fi}}\def\@evenhead{\@oddhead}
\def\@oddfoot{\small\lfoot\ifnum\count0=\startpage\copyright\ \gtp\hfill\else
\agt, Volume \thevolumenumber\ (\thevolumeyear)\hfill\fi}
\def\@evenfoot{\@oddfoot}
\makeatother
\fi
\let\maketitlepage\makeagttitle
\let\makeshorttitle\maketitlepage
\let\maketitle\maketitlepage

\newwrite\gtoutfile
\long\gdef\makeheadfile{  
{\def\\{, }\def\s{ }
\immediate\openout\gtoutfile head.xxx
\immediate\write\gtoutfile{Proxy-for: \ifx\theasciiauthors\relax
\theauthors\else\theasciiauthors\fi\s<\ifx\theasciiemail\relax\theemail\else\theasciiemail\fi>}
\immediate\write\gtoutfile{\noexpand\\}
\immediate\write\gtoutfile{Authors: \ifx\theasciiauthors\relax
\theauthors\else\theasciiauthors\fi}
{\def\\{ }\immediate\write\gtoutfile{Title: \ifx\theasciititle\relax
\thetitle\else\theasciititle\fi}}
\immediate\write\gtoutfile{Subj-class: GT or SG, GR etc}
\immediate\write\gtoutfile{MSC-class: \theprimaryclass\ifx\thesecondaryclass\relax\else, \thesecondaryclass\fi}
\immediate\write\gtoutfile{Journal-ref: Algebr. Geom. Topol. \thevolumenumber\s
(\thevolumeyear) \startpage-\finishpage}
\immediate\write\gtoutfile{Comments: Published by Algebraic and
Geometric Topology at}
\immediate\write\gtoutfile{\s\s\s  http://www.maths.warwick.ac.uk/agt/AGTVol\thevolumenumber/agt-\thevolumenumber-\thepapernumber.abs.html}
\immediate\write\gtoutfile{\noexpand\\}
\immediate\write\gtoutfile{}
\ifx\theasciiabstract\relax
\immediate\write\gtoutfile{\theabstract}\else
\immediate\write\gtoutfile{\theasciiabstract}\fi
\immediate\write\gtoutfile{}
\immediate\write\gtoutfile{\noexpand\\}
\immediate\write\gtoutfile{}
\immediate\closeout\gtoutfile}}  

\def\maketitlepage{\makeagttitle\makeheadfile}
\let\makeshorttitle\maketitlepage
\let\maketitle\maketitlepage

\lognumber{39}
\volumenumber{4}
\volumeyear{2004}
\papernumber{39}
\pagenumbers{893}{934}
\received{17 January 2004} 
\revised{13 September 2004}
\accepted{19 September 2004}
\published{13 October 2004}

\usepackage{graphicx}
\usepackage{amsmath,amssymb}

\let\Bbb\mathbb

\theoremstyle{plain}
\newtheorem*{theorem*}{Theorem}
\newtheorem*{lemma*} {Lemma}
\newtheorem*{corollary*} {Corollary}
\newtheorem*{proposition*} {Proposition}
\newtheorem{theorem}{Theorem}[section]
\newtheorem{lemma}[theorem]{Lemma}
\newtheorem{corollary}[theorem]{Corollary}
\newtheorem{proposition}[theorem]{Proposition}
\newtheorem{question}[theorem]{Question}

\theoremstyle{remark}
\newtheorem*{remark}{Remark}
\newtheorem*{definition}{Definition}

\theoremstyle{definition}

\def \R {\mathbf{R}}
\def \Z {\mathbf{Z}}
\def \C {\mathbf{C}}

\def \L {\mathbf{L}}
\def \G {\mathbf{G}}

\def\eps{\epsilon}

\def\s{\sigma}

\def\Q{\Bbb{Q}}
\def\id{\mbox{id}}

\def\sign{\mbox{sign}}
\def\Z{\Bbb{Z}}
\def\C{\Bbb{C}}

\def\N{\Bbb{N}}
\def\l{\lambda}
\def\part{\partial}

\def\a{\alpha}

\def\tor{\mbox{Tor}}
\def\bp{\begin{pmatrix}}
\def\Arf{\mbox{Arf}}
\def\sm{\setminus}
\def\ep{\end{pmatrix}}
\def\bn{\begin{enumerate}}

\def\en{\end{enumerate}}
\def\ba{\begin{array}}
\def\ea{\end{array}}

\def\L{\Lambda}

\def\s{\sigma}
\def\a{\alpha}
\def\b{\beta}
\def\ti{\tilde}
\def\lk{\mbox{lk}}

\def\fr12{\frac{1}{2}}
\def\diag{\mbox{diag}}

\def\Aut{\mbox{Aut}}
\def\COT{Cochran, Orr and Teichner }
\def\CG{Casson and Gordon }
\def\Oplus{\bigoplus}

\def\ker{\mbox{Ker}}

\def\t{\theta}

\def\hom{\mbox{Hom}}
\def\Ext{\mbox{Ext}}

\def\gcd{\mbox{gcd}}
\def\arf{\mbox{Arf}}

\def\Gal{\mbox{Gal}}

\def\G{\Gamma}
\def\lteta{\eta^{(2)}}
\def\lp{linking pairing }

\def\blp{Blanchfield pairing }
\def\lpp{linking pairing}

\def\blpp{Blanchfield pairing}
\def\ds{doubly slice }
\def\zeros{zeros }

\def\varnothing{\emptyset}

\def\modd{\, mod \, }

\begin{document}

\title[Eta invariants as sliceness obstructions]{Eta invariants as sliceness obstructions and\\their relation to Casson-Gordon invariants}
\author{Stefan Friedl}

\address{Department of Mathematics, Rice University, Houston, TX 77005, USA}
\email{friedl@rice.edu}
\urladdr{http://math.rice.edu/~friedl/}
\asciiurl{http://math.rice.edu/ friedl/}

\begin{abstract}
We give a useful classification of the metabelian unitary
representations of $\pi_1(M_K)$, where $M_K$ is the result of
zero-surgery along a knot $K\subset S^3$.  We show that certain eta
invariants associated to metabelian representations $\pi_1(M_K)\to
U(k)$ vanish for slice knots and that even more eta invariants vanish
for ribbon knots and doubly slice knots.  We show that our vanishing
results contain the Casson--Gordon sliceness obstruction.  In many
cases eta invariants can be easily computed for satellite knots.  We
use this to study the relation between the eta invariant sliceness
obstruction, the eta-invariant ribbonness obstruction, and the
$L^2$--eta invariant sliceness obstruction recently introduced by
Cochran, Orr and Teichner.  In particular we give an example of a knot
which has zero eta invariant and zero metabelian $L^2$--eta invariant
sliceness obstruction but which is not ribbon.
\end{abstract}

\asciiabstract{%
We give a useful classification of the metabelian unitary
representations of pi_1(M_K), where M_K is the result of
zero-surgery along a knot K in S^3.  We show that certain eta
invariants associated to metabelian representations pi_1(M_K) -->
U(k) vanish for slice knots and that even more eta invariants vanish
for ribbon knots and doubly slice knots.  We show that our vanishing
results contain the Casson-Gordon sliceness obstruction.  In many
cases eta invariants can be easily computed for satellite knots.  We
use this to study the relation between the eta invariant sliceness
obstruction, the eta-invariant ribbonness obstruction, and the
L^2-eta invariant sliceness obstruction recently introduced by
Cochran, Orr and Teichner.  In particular we give an example of a knot
which has zero eta invariant and zero metabelian L^2-eta invariant
sliceness obstruction but which is not ribbon.}

\asciikeywords{Knot concordance, Casson-Gordon invariants, Eta invariant}
\maketitle

\section{Introduction}
\subsection{A quick trip through knot concordance theory}
A knot $K\subset S^{n+2}$ is a smooth oriented submanifold homeomorphic to $S^n$.
A knot $K$  is called   slice if it
bounds a smooth $(n+1)$-disk in $D^{n+3}$.
Isotopy classes of knots form a semigroup under connected sum.
The quotient of this semigroup by the subsemigroup of slice knots turns out to be a group,
called the knot concordance group. For example, the existence of inverses
follows by noting that the connected sum of $K$ and $-rK$ bounds a disk,
where $-rK$ denotes the knot obtained by reflecting $K$ through a disjoint hypersphere
and reversing the orientation.
It is a natural goal to attempt to understand this group and to
find complete  invariants for
detecting when a knot is slice.

Knot concordance in the high-dimensional case is well understood.
For even $n$ Kervaire \cite{K65} showed that every knot $K\subset S^{n+2}$ is slice,
and for odd $n\geq 3$ Levine \cite{L69}
showed that $K\subset S^{n+2}$ is slice if and only if
$K$ is algebraically slice (cf.\ section \ref{sectionslice} for a definition). In  this way, the task of detecting slice knots in higher dimensions was reduced to an algebraic problem which is well-understood \cite{L69b}.

The classical case $n=1$ turns out to be a much more difficult problem.
\CG \cite{CG78}, \cite{CG86} defined certain sliceness obstructions
(cf.\ section \ref{sectioncg1})
and used these to give examples of knots in
$S^3$ which are algebraically slice but not geometrically slice.

 Letsche \cite{L00} introduced two new approaches to finding obstructions
 to the sliceness of a knot.
 One approach used the concept of a universal group to find representations that
 extend over the complement of a
 slice disk, the other used eta invariants of metabelian representations of knot
 complements to give sliceness obstructions.
 The idea of universal groups was taken much further in a ground breaking paper
  by Cochran, Orr and Teichner \cite{COT03} (cf.\ section  \ref{sectioncot}).
The goal of this paper is to build on Letsche's second approach.

\subsection{Summary of results}
Given a closed smooth three--manifold $M$ and given a unitary representation
$\a \co \pi_1(M)\to U(k)$
Atiyah--Patodi--Singer
\cite{APS75} defined $\eta(M,\a)\in \R$, called the eta invariant of $(M,\a)$. This invariant is closely related to signatures,
and therefore well--suited for studying cobordism problems.
For a knot $K$
we study the eta invariants associated to the closed manifold $M_K$, the result
of zero--framed surgery along $K\subset S^3$.

In proposition \ref{lemma1} we give a complete classification
of irreducible, unitary, metabelian (cf.\ section \ref{sectionmetab}) representations
of $\pi_1(M_K)$. In theorem \ref{mainthm} we show that for a slice knot the eta invariant vanishes for certain irreducible
metabelian representations of prime power dimensions.

In section \ref{sectioncg1} we recall the Casson--Gordon sliceness
obstruction theorem.
We show in theorem \ref{thmetacg} that for a given knot $K$
the Casson--Gordon sliceness obstruction vanishes
if and only if $K$ satisfies the vanishing conclusion of theorem \ref{mainthm}.

Despite this equivalence of obstructions
the eta invariant approach has several advantages over the Casson-Gordon approach. For example
we show in  theorem \ref{mainthm2} that eta invariants vanish
for tensor products of certain irreducible prime power dimensional representations.
This gives sliceness obstructions, which are   potentially stronger than the Casson--Gordon obstruction.

Recently \COT \cite{COT03}, \cite{COT04}
defined the notion of
$(n)$--solvability for a knot,
$n\in \frac{1}{2}\N$, which has in particular the following properties.
\bn
\item  Slice knots are $(n)$--solvable for all $n$.
\item $(n)$--solvable knots are $(m)$--solvable for all $m\leq n$.
\item A knot is $(0.5)$--solvable if and only if it is algebraically slice.
\item $(1.5)$--solvable knots  have zero Casson--Gordon obstruction.
\en

Given a homomorphism $\varphi\co\pi_1(M)\to G$ to a group $G$,
Cheeger and Gromov defined the $L^2$--eta invariant $\eta^{(2)}(M,\varphi)\in \R$
which is closely related to Atiyah's $L^2$--signature.
 \COT used $L^2$--eta invariants
to find examples of knots which are $(2.0)$--solvable,
which in particular have zero Casson--Gordon obstructions, but which are not $(2.5)$--solvable. Using similar ideas
Taehee  Kim
\cite{K02} found examples of knots which are $(1.0)$--solvable and  have zero Casson--Gordon obstructions, but which
are not
$(1.5)$--solvable. A quick summary of this theory is given in section \ref{sectioncot}.

In section  \ref{sectionex} we give examples that
show that $L^2$--eta invariants are not complete (non--torsion)
invariants for a knot to be $(0.5)$--solvable respectively
$(1.5)$--solvable.
In our examples  we use a satellite
construction to get knots whose eta invariants can be computed explicitly by methods introduced
by Litherland \cite{L84}.

The systematic study of eta invariants corresponding to metabelian representations
also allows us to find
more refined obstructions to a knot being ribbon
(theorem \ref{mainribbonthm}).
It is not known whether all slice knots satisfy the ribbon obstruction theorem.
The fact that the ribbon obstruction is apparently stronger than the corresponding
sliceness obstruction is particularly interesting as this could
potentially provide a way
to disprove the {\em ribbon conjecture } (that every slice knot is ribbon).
 In proposition \ref{propex3} we give examples of knots
for which all abelian and metabelian sliceness obstructions vanish,
 but which do not satisfy the
condition for  theorem \ref{mainribbonthm}, i.e.\
which are not ribbon.
It is not known whether these examples are slice or not.

We also apply our methods to the study of doubly slice knots
(theorem \ref{maindoublyslicethm}). It is a well--known
fact that the doubly sliceness condition is  much stronger than the
sliceness condition  (cf.\ \cite{S71}, \cite{K03}, \cite{K04}).
We point out the interesting fact that doubly slice knots satisfy the ribbon obstruction
theorem, which
suggests that doubly slice knots have a
`higher chance' of being ribbon than ordinary slice knots.

We conclude this paper with a discussion of Gilmer's
\cite{G83}\cite{G93} and Letsche's \cite{L00} obstruction theorems
in sections  \ref{sectiongilmer} and \ref{sectionletsche}.
During our work on eta invariants as sliceness obstructions
we found that both theorems  have gaps in their proofs. We explain
where the problem lies and we show that Gilmer's and Letsche's
approach still give ribbonness obstruction.
These ribbon obstructions can be shown to be contained in theorem \ref{mainribbonthm}.  \\

This paper is essentially the author's thesis.
I would like to express my deepest gratitude towards my advisor Jerry Levine for
his patience,  help and encouragement.

\section{Basic knot theory  and linking pairings}

Throughout this paper we will always work in the smooth category. In the classical dimension the
theory of knots in the smooth
category is equivalent to the theory in the locally flat category.
All homology groups are furthermore taken over $\Z$, unless otherwise indicated.

\subsection{Basic knot theory} \label{sectionbasics}
By a {\em knot} $K$ we understand an oriented submanifold of $S^3$
diffeomorphic to $S^1$. A oriented surface $F \subset S^3$
with $\partial(F)=K$ will be called a {\em Seifert surface }for $K$.
Note that a Seifert surface inherits an orientation from $K$,
in particular the map $H_1(F)\to H_1(S^3\sm F), a\mapsto a_+$
induced by pushing into the positive normal direction is well-defined.
The pairing
\[ \ba{rcl} H_1(F) \times H_1(F) &\to&\Z \\
   (a,b) & \mapsto& \lk(a,b_+) \ea \]
is called the {\em Seifert pairing} of $F$. Any matrix $A$ representing such a
pairing
is called a {\em Seifert matrix } for $K$.
By  \cite[p.\ 393]{M65}, \cite{L70} two Seifert matrices for a given knot $K$ are  S--equivalent.
In particular the Alexander polynomial $\Delta_K(t):=\det(At-A^t) \in \Z[t,t^{-1}]$
is a well-defined invariant of $K$, i.e.\ independent of the choice of $A$, up to multiplication by units $\pm t^l$.

Let $N(K)$ be a solid torus neighborhood of $K$. A meridian $\mu$ of $K$ is a non-separating simple  closed curve in
$\partial(N(K))$ that bounds a disc in $N(K)$.
The notion of meridian  is well-defined up to homotopy in $S^3\sm K$.
Alexander duality shows that there exists an isomorphism $\eps\co H_1(S^3 \sm K)\to \Z$
 such that $\eps(\mu)=1$.

We will mostly study invariants of $M_K$, the result of zero-framed-surgery along $K$.
This manifold
has the advantage over $S^3 \sm K$ that it is
a closed manifold associated to $K$.
Note that we have isomorphisms $H_1(M_K)\xleftarrow{\cong}H_1(S^3\sm K) \xrightarrow{\eps} \Z$, we will denote this composition of maps
by $\eps$ as well.

\subsection{Slice knots} \label{sectionslice}
 
We begin with a few definitions.
We say that two knots $K_1,K_2$ are {\em smoothly concordant }if there exists a smooth
submanifold $V \subset S^3 \times [0,1]$ such that $V \cong S^1 \times [0,1]$
and such that $\partial(V)=K_0 \times 0 \cup -K_1 \times 1$.
A knot $K \subset  S^3$ is called {\em smoothly slice} if it is concordant
to the unknot. Equivalently, a knot is smoothly slice if there exists
 a smooth disk $D \subset D^4$ such that
$\partial(D)=K$.

A knot $K$ is called {\em topologically slice} if $K$ bounds a locally flat disk $D\subset D^4$. Locally flat means that one can find an embedding
$D\times D^2\subset D^4$ which extends the embedding of the disk.

We say that the Seifert pairing on $H_1(F)$ is {\em }metabolic if there
exists a subspace $H$ of half-rank such that the Seifert pairing vanishes on $H$.
This is the case if and only if $K$ has a metabolic Seifert matrix $A$, i.e.\ if $A$ is of the
 form
\[A= \bp 0 & B \\ C & D \ep, \]
where $0,B,C,D$ are square matrices.
If the Seifert pairing of a knot is metabolic then we say that $K$ is algebraically slice.

By Levine \cite{L69} we now get the following implications
\[ \mbox{smoothly slice} \Rightarrow \mbox{topologically slice} \Rightarrow \mbox{algebraically slice}. \]
The reverse implications do not hold. \CG \cite{CG86} gave examples of
algebraically slice knots which are not topologically slice.
Freedman \cite{FQ90} showed that any knot with trivial Alexander
polynomial is topologically slice, whereas Gompf \cite{G86} gave examples
of knots with trivial Alexander polynomial that are not smoothly slice.

In the following we will work in the smooth category. In particular by slice we will mean smoothly slice.

 For a slice disk $D$ we define  $ N_D:=\overline{D^4\sm N(D)}$, where $N(D)$ is a tubular neighborhood of $D$ in $D^4$.
We summarize a couple of well-known facts about $N_D$.
 \bn
\item $\partial(N_D)=M_K$,
\item $H_*(N_D)=H_*(S^1)$,
\item the inclusion map $H_1(M_K)\to H_1(N_D)$ is an isomorphism.
\en
We will denote the induced map $H_1(N_D)\xleftarrow{\cong} H_1(M_K)\xrightarrow{\eps} \Z$ by $\eps$ as well.

\subsection{Universal abelian cover of $M_K$ and the \blp}  \label{sectionabelcov}

\label{sectioninfinitecover}
Let $K$ be a knot.
Denote
the infinite cyclic cover of $M_K$
corresponding to $\eps\co  H_1(M_K)\to \Z$
by  $\ti{M}_K$. Then $\Z=\langle t \rangle$  
acts on   $\ti{M}_K$,  
therefore  
 $H_1(\ti{M}_K)$ carries
a $\L:=\Z[t,t^{-1}]$-module structure.
Clearly $H_1(\ti{M}_K)$ and the twisted homology module 
$H_1(M_K,\L)$ are canonically isomorphic as $\L$--modules.

In the following we will give $\L$ the involution induced by $\bar{t}=t^{-1}$.
Let $S:=\{ f\in \L | f(1)=1\}$. The $\L$-module $H_1(M_K,\L)$ is $S$-torsion
since for example the Alexander
polynomial $\Delta_K(t)$ lies in $S$ and $\Delta_K(t)$ annihilates $H_1(M_K,\L)$.
Blanchfield \cite{B57} introduced the pairing  
\[ \ba{rcl} \l_{Bl}\co H_1(M_K,\L) \times H_1(M_K,\L) & \to &S^{-1}\L/\L  \\
(a,b) &\mapsto & \frac{1}{p(t)} \sum_{i=-\infty}^{\infty} (a\cdot t^ic)t^{-i}
\ea  \]
where $c\in C_2(M_K,\L)$ such that $\partial(c)=p(t)b$ for some $p(t) \in S$. It is well-defined,
 non-singular and hermitian over $\L$.
For any $\L$-submodule $P\subset H_1(M_K,\L)$ define
\[ P^{\perp}:=\{ v\in H_1(M_K,\L) | \l_{Bl}(v,w)=0 \mbox{ for all } w \in P\}.
\]
If $P \subset H_1(M,\L)$ is such that $P=P^{\perp}$ then we say
that $P$ is a metabolizer for $\l_{Bl}$. If $\l_{Bl}$ has a metabolizer we say that
$\l_{Bl}$ is metabolic.

For a $\Z$--module $A$ denote by  $TA$ the $\Z$--torsion
submodule of $A$ and let $FA:=A/TA$.

\begin{theorem}[{\cite{K75} \cite[prop.\ 2.8]{L00}}]
\label{keartonthm}    \label{thmblanchfield1}   \label{firstpropmetabl}
$\phantom{99}$
\bn
\item
If $K$ is slice and $D$ any  slice disk, then  $$P:=\ker\{H_1(M_K,\L)\to FH_1(N_D,\L)\}$$ is a metabolizer
for the Blanchfield pairing.
\item
The \blp for a knot $K$ is
metabolic if and only if
$K$ is algebraically slice.
\en
\end{theorem}

It is an open problem whether every slice knot
$K$ admits a slice disk $D$ such that
the submodule $\ker\{H_1(M_K,\L)\to H_1(N_D,\L)\}$ is a metabolizer for
the Blanchfield pairing.
This is an important question since it is closely linked to the ribbon conjecture
(cf.\ \ref{sectionrib}).

\subsection{Finite cyclic covers and linking pairings}
\label{appendixcycliccovers}
Let $K$ be a knot, $k$ some number.
Denote by
$M_k$ the $k$-fold cover of $M_K$ corresponding to
$\pi_1(M_K)\xrightarrow{\eps} \Z \to \Z/k$, and denote by
$L_k$ the $k$-fold cover of $S^3$ branched along $K \subset S^3$.
Note that $H_1(L_k)$ and $H_1(M_k)$ also have $\L$--module structures.

\begin{lemma} \label{lemmah1mkh1lk}  \label{propofcycliccov}
$$ |H_1(L_k)| = \left| \prod_{j=1}^k \Delta_K(e^{2 \pi i j/k})\right|,
\leqno{\phantom{9}\rm(1)}$$
 where $0$ on the right hand side means that $H_1(L_k)$ is infinite.
\leftskip 28.5pt

\bn
\item[\rm(2)]
There exist natural isomorphisms  
\[ \ba{rcl}  
H_1(L_k) &=& H_1(M_K,\L)/(t^k-1),\\
 H_1(M_k) &=& H_1(L_k)\oplus \Z =H_1(M_K,\L)/(t^k-1) \oplus \Z.
\ea \]
\en
\end{lemma}

\begin{proof}
The first part is shown in  \cite[p.\ 17]{G77}. For the second part consider
  the long exact sequence
\[ \dots \to H_i(\ti{M}_K)\xrightarrow{t^k-1} H_i(\ti{M}_K)\to H_i(M_k)\to \dots.\]
It follows immediately that $H_1(M_k) =H_1(M_K,\L)/(t^k-1) \oplus \Z$.
The other isomorphism is clear.
\end{proof}

Now let $k$ be any integer  such that $H_1(L_k)$ is finite, then the map
\[ \ba{rcl} \l_{L,k}=\l_L\co H_1(L_k) \times H_1(L_k) &\to &\Q/\Z \\
     (a,b) &\mapsto & \frac{1}{n} a\cdot c \mod \Z, \ea \]
where $c \in C_2(L_k)$ such that $\partial(c)=nb$,  
defines a symmetric, non-singular pairing which is called the \lp  
 of $H_1(L_k)=TH_1(M_k)$.
For any $\Z$-submodule $P_k\subset TH_1(M_k)$ define  
\[ P_k^{\perp}:=\{ v\in TH_1(M_k) | \l_{L}(v,w)=0 \mbox{ for all } w \in P_k\}. \]
If $P_k$ is a $\L$--submodule with $P_k=P_k^{\perp}$ then we say that $P_k$ is a metabolizer
for
the linking pairing $\l_L$. If $\l_L$ has a metabolizer then we say
that $\l_L$ is metabolic.
Note that if $P_k$ is a metabolizer for the linking pairing, then $|P_k|^2=|H_1(L_k)|$.

\subsection{Homology of prime-power covers and the \lp}

The following corollary shows that prime power covers of manifolds behave in a more `controlled' way.

\begin{lemma}\label{h1mk}
Let $Y$ be a  manifold such that $H_*(Y)=H_*(S^1)$
and let $k=p^s$ where $p$ is a prime number. Write $Y_k$ for the $k$--fold cover of $Y$ corresponding to $H_1(Y)=H_1(S^1)=\Z \to \Z/k$.
Then $H_*(Y_k)=H_*(Y) \oplus \mbox{torsion}$.
\end{lemma}

The following proof is modelled after  \cite[p.\ 184]{CG86}.

\begin{proof}
For any $n \in \N \cup \{\infty\}$ we can give $H_*(Y_n)$ a $\L$-structure.
Since $H_2(Y,\Z/p)=0$ we get the following exact sequence
\[  0 \to H_1(Y_{\infty},\Z/p) \xrightarrow{t-1}
  H_1(Y_{\infty},\Z/p) \to H_1(Y,\Z/p) \to H_0(Y_{\infty},\Z/p) \to 0.
\]
Since $ H_1(Y,\Z/p) \to H_0(Y_{\infty},\Z/p)$ is an isomorphism it follows
from the sequence that  the map $H_1(Y_{\infty},\Z/p) \xrightarrow{t-1}
H_1(Y_{\infty},\Z/p)$ is an isomorphism. Since $H_i(Y)=0$ for $i>1$ we also
get that
$H_i(Y_{\infty},\Z/p) \xrightarrow{t-1}     H_i(Y_{\infty},\Z/p)$ is an isomorphism for $i>1$.
Over $\Z/p$ we  get  $(t^k-1)=(t^{p^s}-1)=(t-1)^{p^s}$, hence multiplication by $(t^k-1)$ is an
automorphism of $H_1(Y_{\infty},\Z/p)$ as well.  Consider the long exact sequence  
\[ 
\dots \to H_1(Y_{\infty},\Z/p) \xrightarrow{t^k-1} 
  H_1(Y_{\infty},\Z/p) \to H_1(Y_k,\Z/p) \to H_0(Y_{\infty},\Z/p) \to 0.
\]
It follows that $H_1(Y_k,\Z/p)=H_0(Y_{\infty},\Z/p)=\Z/p$.
Similarly we can show that $H_i(Y_k,\Z/p)=0$ for $i>1$.
The lemma now follows from the universal coefficient theorem.  
\end{proof}

\begin{corollary} \label{corh1mk}  \label{proplinkingformmetabolic}  
Let $k$ be a prime power.
$$ \ba{rcl} H_1(M_k)&=&\Z \oplus TH_1(M_k).\\
TH_1(M_k)&=&H_1(L_k)=H_1(M_K,\L)/(t^k-1).  \ea \leqno{\phantom{9}\rm(1)}$$
\leftskip28.5pt
\bn
\item[\rm(2)] Let $D$ be a slice disk, denote the $k$-fold cover of $N_D$ by $N_k$, then
$H_1(N_k)=\Z \oplus TH_1(N_k)=\Z \oplus H_1(N,\L)/(t^k-1)$.
\item[\rm(3)] $Q_k:=\ker\{  TH_1(M_k) \to TH_1(N_k)\}$ is a metabolizer for $\l_L$.
\en
\end{corollary}  

The first two statements are immediate corollaries. The third statement is well--known (cf.\ \cite{G77} or \cite{F03d}). It is important though
that $k$ is a prime power, since it is crucial in the proof that $H_1(N_k)=\Z \oplus TH_1(N_k)$.

\begin{remark}
Let $K$ be a slice knot, $D$ a slice disk. Then $P:=\ker\{H_1(M_K,\L)$\break
$\to FH_1(N_D,\L)\}$ is a metabolizer for
the Blanchfield pairing and $Q_k:=$\break
$\ker\{  TH_1(M_k) \to TH_1(N_k)\}$ is a metabolizer for $\l_{L,k}$.
Furthermore $P_k:=\pi(P)\subset H_1(M_K,\L)/(t^k-1)$ can be shown to be a metabolizer  for $\l_{L,k}$.
It is an open problem whether $P_k=Q_k$. If yes, then this would show that all the metabolizers $Q_k$
can be lifted to a metabolizer for the Blanchfield pairing, which would resolve the problems which appeared
in \cite{G93} and \cite{L00} (cf.\ sections \ref{sectiongilmer} and \ref{sectionletsche}).
\end{remark}

\section{Introduction to eta-invariants and first application to knots}

\subsection{Eta invariants}

Let $C$ be a complex, hermitian matrix, i.e.\ $C=\bar{C}^t$, then   the signature $\sign(C)$ is defined
as the number of positive eigenvalues
of $C$ minus the number of negative eigenvalues.
The following is an easy exercise.  
\begin{lemma}  \label{lemmasginzero}
If $C$ is hermitian then $\sign(PC\bar{P}^t)=\sign(C)$ for any 
$P$ with $\det(P)\ne 0$. If furthermore $\det(C)\ne 0$ and $C$ is of the form  
\[C= \bp 0 & B \\ \bar{B}^t & D \ep, \]
where $B$ is a square matrix,
then $\sign(C)=0$.  
\end{lemma}

Let $M^{3}$ be a closed  manifold 
and $\a\co\pi_1(M) \to U(k)$ a representation. Atiyah,
Patodi, Singer \cite{APS75} associated to $(M,\a)$ a number $\eta(M,\a)$
called the (reduced) eta invariant of
$(M,\a)$.

The main theorem to compute the eta invariant is the following.

\begin{theorem}[Atiyah--Patodi--Singer index theorem \cite{APS75}] \label{apsthm}
 If there exists a manifold $W^{4}$ and a representation $\b\co\pi_1(W)\to U(k))$
 such that  $\partial(W,\b)=r(M^3,\a)$
for some $r\in \N$, then
\[ \eta(M,\a) =\frac{1}{r}(\sign_{\b}(W)-k\, \sign(W)), \]
where $\sign_{\b}(W)$ denotes the   signature of the twisted intersection pairing
on $H_2^{\b}(W)$.
\end{theorem}

\subsection{Application of eta invariants to knots}
\label{sectionfirstappl}
For a knot $K$
we will study the eta invariants associated to the closed manifold $M_K$.
In the context of knot theory they were studied by Levine
\cite{L94} who used them
to find links which are not concordant to boundary links.
Letsche \cite{L00} used eta invariants to study knot concordance in dimension three.

\begin{definition}
For a group $G$ the derived series is defined by
$G^{(0)}:=G$ and inductively
$G^{(i+1)}:=[G^{(i)},G^{(i)}]$ for $i>0$.
\end{definition}

The inclusion map $S^3 \sm K \to M_K$ defines a homomorphism
$\pi_1(S^3 \sm K) \to \pi_1(M_K)$,
the kernel is generated by the longitude of $K$ which lies in $\pi_1(S^3\sm K)^{(2)}$.
In particular
$\pi_1(S^3\sm K)/\pi_1(S^3\sm K)^{(2)}\to\pi_1(M_K)/\pi_1(M_K)^{(2)}$
is an isomorphism.

Let $M_O$ be the zero-framed surgery on the trivial knot.
Then $M_O=S^1\times S^2$ which bounds
$S^1\times D^3$, which is homotopic to a 1-complex. This proves that the unknot has vanishing eta invariant for 
any unitary representation since any representation of $M_O$ extends over $S^1\times D^3$.

\subsection{$U(1)$-representations} \label{sectionu1rep}

 Let $K$ be a knot, $A$ a Seifert matrix, then
we define the signature function $\s(K)\co S^1 \to \Z$ of $K$
as follows (cf.\ \cite[p.\ 242]{L69})
\[ \ba{rcl} \s_z(K) :=\sign(A(1-z)+A^t(1-\bar{z})).
      \ea \]
It is easy to see that this is independent of the choice of $A$.

\begin{proposition} \label{propu1sign}
Let $K$ be a knot, $\mu$ a meridian and
let $\a \co \pi_1(M_K) \to U(1)$ be a representation.
\bn
\item Let $z:=\a(\mu)$, then
\[ \eta(M_K,\a)=\s_z(K). \]
 \item The function $z\mapsto \s_z(K)$ is
 locally constant outside of the zero set of $\Delta_K(t)$.
\item If $K$ is
algebraically slice and $z \in S^1$ such that $\Delta_K(z) \ne 0$, then
$\s_{z}(K)=0$. In particular if $z$ is a prime power root of unity, then
$\s_z(K)=0$.
\item Given $z$ a prime power root of unity, $K\mapsto \s_z(K)$ defines a homomorphism from the knot concordance group to $\Z$.
\en
\end{proposition}

\proof
\bn
\item See \cite{L84}.
\item The function $z \mapsto \s_z(K)$ is continuous on
\[ \{ z\in S^1 | A(1-z)+A^t(1-\bar{z}) \mbox{ is non--singular} \}.\]
an easy argument shows that this set equals
\[ \{1\} \cup \{ z\in S^1 | \Delta_K(z)=\det(Az-A^t)\ne 0 \}.\]
Levine \cite{L69} showed that the signature function is continuous at $z=1$.
\item If $z \in S^1$ such that $\Delta_K(z) \ne 0$ then $A(1-z)+A^t(1-\bar{z})$ is non--singular
and the first part follows from lemma  \ref{lemmasginzero}.  The second part follows from
$\Delta_K(1)=1$ and the well--known fact that   $\Phi_{p^r}(1)=p$ for a prime $p$,
where $\Phi_{p^r}(t)$ denotes the minimal polynomial of
a primitive root of unity of order $p^r$.
\item Follows immediately from (3) and the additivity of the twisted signature function.\qed
\en

Let $K$ be a knot, $m\in \N$, then we can form $mK$ by iterated connected sum.
We say  $K$ is (algebraically) torsion if $mK$ is (algebraically) slice for some $m$.
It follows from proposition  \ref{propu1sign} (4) that $\s_z(K)=0$ for any algebraically torsion knot $K$ and $z$
a prime power root of unity.

Levine \cite{L69} showed that for $n>1$ a knot $K\subset S^{2n+1}$ is slice
if and only if $K$ is algebraically slice. Levine \cite{L69b} and  Matumuto \cite{M77}
showed that
$K$ is algebraically torsion if and only if $\eta(M_K,\a)=0$
for all $\a\co\pi_1(M_K)\to U(1)$ of prime power order.
In particular, in the case $n > 1$ the $U(1)$-eta invariant detects any
non--torsion knot.

In the classical case $n=1$ Casson and Gordon \cite{CG86} found
an example of a non--slice knot which is algebraically slice.
Jian \cite{J81} showed using the Casson--Gordon invariant that there
are knots which are algebraically slice and which are not torsion in the knot
concordance group.

The goal of this paper is to study to which degree
certain non-abelian eta invariants can detect algebraically slice knots which are not
topologically slice.

\section{Metabelian eta invariants and the main sliceness obstruction theorem}
\label{sectionmetab}
A group $G$ is called metabelian if $G^{(2)}=\{e\}$, a
 representation $\varphi\co\pi_1(M)\to U(k)$ is called
metabelian if it factors through $\pi_1(M)/\pi_1(M)^{(2)}$.
Metabelian eta invariants in the context of knot concordance were  studied by Letsche \cite{L00}.
Our approach is influenced by his work.

For a knot $K$ let $\pi:=\pi_1(M_K)$, then consider
\[ 1\to \pi^{(1)}/\pi^{(2)}\to \pi/\pi^{(2)}\to\pi/\pi^{(1)}\to 1. \]
Note that $\pi_1(\ti{M}_K)=\pi_1(M_K)^{(1)}$ where $\ti{M}_K$ denotes
the infinite cyclic cover of $M_K$, hence
$H_1(M_K,\L)=H_1(\ti{M}_K) \cong \pi_1(M_K)^{(1)}/\pi_1(M_K)^{(2)}$.
Since $\pi/\pi^{(1)}=H_1(M_K)=\Z$ the above sequence splits and we get  isomorphisms
\[  \pi/\pi^{(2)}
\cong  \pi/\pi^{(1)} \ltimes \pi^{(1)}/\pi^{(2)}
\cong   \Z \ltimes H_1(M_K,\L),  \]
where $n \in \Z$ acts by multiplication
by $t^n$. This shows that metabelian representations of $\pi_1(M_K)$ correspond to representations of $\Z \ltimes H_1(M_K,\L)$.

\subsection{Metabelian representations of $\pi_1(M_K)$} \label{sectionclassrep}
\label{sectionmetabrep}
For a group $G$ denote by
$R_k^{irr}(G)$ (resp. $R_k^{irr,met}(G)$) the
set of conjugacy classes of irreducible, $k$-dimensional,  unitary
(metabelian)
representations of
$G$. Note that the eta invariant of a manifold only depends on the conjugacy class of a unitary representation.
Recall that for a knot $K$ we can identify
\[ R_k^{irr,met}(\pi_1(M_K)) = R_k^{irr}(\Z \ltimes H_1(M_K,\L)). \]
The following proposition gives a very useful classification of all irreducible metabelian unitary representations
of $\pi_1(M_K)$ (cf.\ \cite{L95} for a different approach).

\begin{proposition} \label{lemma1}\label{lemma2}
Let $H$ be a $\L$--module.
Let $z \in S^1$ and $\chi\co H \to H/(t^k-1) \to S^1$ a character, then
 \[ \ba{rcl} \a_{(z,\chi)} \co\Z \ltimes H& \to &U(k) \\
    (n,h) &\mapsto &
z^n \bp 0& \dots &0&1 \\ 1&\dots &0&0 \\
\vdots &\ddots &&\vdots \\
     0&\dots &1&0 \ep^n \bp \chi(h) &0&\dots &0 \\
 0&\chi(th) &\dots &0 \\
\vdots &&\ddots &\vdots \\ 0&0&\dots &\chi(t^{k-1}h) \ep \ea \]
defines a  representation.
If $\chi$ does not factor through $H/(t^l-1)$ for some $l<k$, then $\a_{(z,\chi)}$ is irreducible.

Conversely, any  $[\a] \in{R}_k^{irr}(\Z \ltimes H)$ has a representative
$\a_{(z,\chi)}$ with $z,\chi$ as above.  
\end{proposition}

\begin{proof}
It is easy to check that $\a_{(z,\chi)}$ is well-defined. Now
let $[\a] \in {R}_k^{irr}(\Z \ltimes H)$.
Denote by $\chi_1,\dots,\chi_l\co H\to S^1$  the different weights of $\a\co0\times H\to U(k)$.
Since $H$ is an abelian group we can write
$\C^k=\oplus_{i=1}^l V_{\chi_i}$
where $V_{\chi_i}:=\{ v\in \C^k | \a(0,h)(v)=\chi_i(h)v \mbox{ for all }h\}$ is the weight space
corresponding to $\chi_i$.

Recall that the group structure of $\Z \ltimes H$ is given by
\[ (n,h)(m,k)=(n+m,t^mh+k). \]
In particular   for all $v \in H$
\[ (j,0)(0,t^jh)=(j,t^jh)=(0,h)(j,0), \]
therefore for $A:=\a(1,0)$ we get
\[ A^j\a(0,t^jh)=\a(j,0)\a(0,t^jh)
 =\a(j,t^jh)=\a(0,h)\a(j,0)=\a(0,h)A^j. \]
This shows that $\a(0,t^jh)=A^{-j}\a(0,h)A^j$.
Now let $v\in V_{\chi(h)}$, then
\[ \a(0,h)Av=A\a(0,th)A^{-1}Av=A\a(0,th)v=A\chi(th)v=\chi(th)Av, \]
i.e.\ $\a(1,0)\co V_{\chi_i(h)}\to V_{\chi_i(th)}$. Since $\a$ is irreducible it follows that, after
reordering,
$\chi_i(v)=\chi_1(t^iv)$  for all $i=1,\dots,l$. Note that $A^j$ induces isomorphisms between the weight spaces
$V_{\chi_i}$ and
that
$A^l\co V_{\chi_1}\to V_{\chi_1}$ is a unitary transformation. In particular it has an eigenvector $v$, hence
$\C v \oplus \C Av\oplus\dots \C A^{l-1}v$ spans an $\a$-invariant subspace. Since $\a$ is irreducible it follows
that $l=k$ and that each $V_{\chi_i}$ is one-dimensional.

Since $\a$ is a unitary representation we can find a unitary matrix $P$ such that
$P\C e_i=V_i$, in particular, $\a_1:=P^{-1}\a P$ has the following properties.
\bn
\item $\a(0 \times H)=\diag(\chi(h),\chi(th),\dots,\chi(t^{k-1}h))$,
\item for some $z_1,\dots,z_k \in S^1$ 
\[ \a(1,0):=\bp 0& \dots &0&z_k \\ z_1&\dots &0&0 \\
\vdots &\ddots &&\vdots \\
     0&\dots &z_{k-1}&0 \ep \]
\en
Here we denote by $\diag(b_1,\dots,b_k)$ the diagonal matrix with entries
$b_1,\dots,b_k$.

Let $z:=\prod_{i=1}^k z_i$
and let $Q:=\diag(d_1,\dots,d_k)$  where $d_i:=\frac{\prod_{j=1}^{i-1}z_j}{z^{i-1}}$.
Then $\a_2:=Q^{-1}\a_1Q$ has the required properties.
\end{proof}

\subsection{Eta invariants as concordance invariants} \label{sectionspecialvarieties}
We quote some definitions, initially introduced by Levine \cite{L94}.  
Let $G$ be a group, then a $G$-manifold is a pair $(M,\a)$ where $M$ is a compact oriented manifold 
and $\a$ is
a  homomorphism $\a\co\pi_1(M) \to G$  defined up to inner automorphism.

We call two   $G$-manifolds $(M_j,\a_j), j=1,2$,
homology $G$--bordant if there exists a
$G$-manifold
$(N,\b\co\pi_1(N)\to G)$  such that $\partial(N)=M_1 \cup -M_2, H_*(N,M_j)=0$ for $j=1,2$
and, up to inner automorphisms of $G$,
$\beta|\pi_1(M_j)=\a_j$.

We will compare eta invariants for homology $G$--bordant $G$--manifolds.

\begin{definition}
For a $\L$-torsion module $H$ define
$P_k^{irr}(\Z \ltimes H)$ to be the set
of conjugacy classes of representations which are conjugate to
$\a_{(z,\chi)}$ with $z\in S^1$ transcendental and
$\chi\co H/(t^k-1) \to S^1$ factoring through a group of prime power order.
If $W$ is a manifold with $H_1(W)\cong \Z$ then we define
$P_k^{irr,met}(\pi_1(W)):=P_k^{irr}(\Z \ltimes H_1(W,\L))$.
\end{definition}

We need the following theorem, which is  a slight
reformulation of a theorem by Letsche \cite{L00} which in turn is based on work
by Levine \cite{L94}.

\begin{theorem} \label{letschepktheorem}
Let $H$ be a $\L$-torsion module, $G:=\Z \ltimes H$. If $(M_1,\!\a_1),\! 
(M_2,\!\a_2)$ are homology $G$--bordant
3--manifolds and if $\t
\in
P^{irr}_k(G)$, then $\eta(M_1,\t \circ \a_1)=\eta(M_2,\t \circ \a_2)$.
\end{theorem}

\proof
Write $\mathcal{P}_k(\pi_1(M))$ for the set of unitary representations $\a\co\pi_1(M)\to U(k)$ with the following two properties
(cf.\ \cite[p.\ 311]{L00})
\bn
\item $\a$ factors through a non--abelian group 
of the form $\Z \ltimes P$, $P$ a finite $p$--group and $\Z\ltimes P\to \Z$ induces an isomorphism on first homology,
\item there exists $g\in \pi_1(M)$ that generates $\pi_1(M)$ such that all eigenvalues of $\a(g)$ are transcendental.
\en

Letsche  \cite[prop.\ 1.7, cor.\ 3.10, thm.\ 3.11]{L00} showed that the statement holds for all $\t\in \mathcal{P}_k(G)$. 
Clearly if $\a \in \mathcal{P}_k(G)$, then all its conjugates lie in $\mathcal{P}_k(G)$ as well.

It therefore suffices to show that  $\a_{(z,\chi)} \in \mathcal{P}_k(\Z \ltimes H)$
if  $[\a_{(z,\chi)}]\in P_k^{irr}(\Z\ltimes H)$.
Let
$[\a_{(z,\chi)}]\in P_k^{irr}(\Z \ltimes H)$, i.e.\ $z\in S^1$ transcendental and
$\chi\co H/(t^k-1) \to \Z/m\to S^1$
where $m$ is a prime power.
Then the result
follows immediately
from
the observations that
\bn
\item all the eigenvalues of $\a_{(z,\chi)}(1,0)$ are of the form
$ze^{2\pi ij/k}$, in particular all are transcendental,
\item $\a_{(z,\chi)}\co\Z \ltimes H \to U(k)$ factors through
$\Z \ltimes (\Z/m)^k$ and $(\Z/m)^k$ is a group of prime power order,
 where $\Z$ acts on $(\Z/m)^k$ by cyclic
permutation, i.e.\ by $1\cdot (v_1,\dots,v_k):= (v_k,v_1,\dots,v_{k-1})$,
\item $H_1(\Z  \ltimes (\Z/m)^k)\to H_1(\Z)$
is an isomorphism.\qed
\en

Let $K$ be a slice knot with slice disk $D$ and
let $\a_{(z,\chi)}\in R_k^{irr,met}(\pi_1(M_K))$.
Consider the following diagram
\[ \ba{ccccccl} \pi_1(M_K)&\to &\Z\ltimes H_1(M_K,\L)&\to&\Z\ltimes
H_1(M_K,\L)/(t^k-1)&\xrightarrow{\a_{(z,\chi)}}&U(k)\\
\downarrow &&\downarrow &&\downarrow&&\downarrow \\
\pi_1(N_D)&\to &\Z\ltimes H_1(N_D,\L)&\to&\Z\ltimes
H_1(N_D,\L)/(t^k-1)&&U(k)\ea \]
If $\chi$ vanishes on $\ker\{H_1(M_K,\L)/(t^k-1)\to H_1(N_D,\L)/(t^k-1)\}$
then $\chi$ extends to $\chi_N\co H_1(N_D,\L)/(t^k-1)\to S^1$ since $S^1$ is divisible.
Furthermore, if $\chi$ is of prime power order, then $\chi_N$ can be chosen to be of prime power order as well.
Note that $\a_{(z,\chi_N)}\co\pi_1(N_D)\to U(k)$ is an extension of $\a_{(z,\chi)}\co\pi_1(M_K)\to U(k)$.
This proves the following.
\begin{lemma}\label{lemmaextntond}
\bn
\item $\a_{(z,\chi)}$ extends to a metabelian representation of $\pi_1(N_D)$ if and only if
$\chi$ vanishes on $\ker\{H_1(M_K,\L)/(t^k-1)\to H_1(N_D,\L)/(t^k-1)\}$.
\item If $\a_{(z,\chi)}{\in} P_k^{irr}(\pi_1(M_K))$  extends to a metabelian representation of $\pi_1(N_D)$
then $\a_{(z,\chi)}$ extends to a representation in $P_k^{irr}(\pi_1(N_D))$.
\en
\end{lemma}

\begin{theorem} \label{letschepkcor}
Let $K$ be a slice knot and $D$ a slice disk. If $\a$  extends to a metabelian representation of
$\pi_1(N_D)$ and if
$\a \in P_k^{irr}(\pi_1(M_K))$, then $\eta(M_K,\a)=0$.
\end{theorem}

\begin{proof}
By lemma  \ref{lemmaextntond} we can find an extension $\b\in P_k^{irr}(\pi_1(N_D))$.
We can decompose $N_D$ as $N_D=W^4 \cup_{M_O} S^1 \times D^3$ where
$M_O=S^1 \times S^2$ is the zero-framed surgery along the
trivial knot in $S^3$ and $W$ is a homology $\Z$--bordism
between $M_K$ and $M_O$.
 The statement now follows from theorem \ref{letschepktheorem}
and the fact that the unknot has zero eta invariants,
since $(W,\id)$ is a homology $\Z \ltimes H_1(W,\L)$--bordism between
$(M_K,i_*)$ and $(M_O,i_*)$.
\end{proof}

\subsection{Main sliceness obstruction theorem}
\begin{theorem} \label{mainthm}
Let $K$ be a slice knot, $k$ a prime power. Then there
exists a metabolizer $P_k \subset TH_1(M_k)$ for the \lpp, such that
for any
  representation $\a\co\pi_1(M_K)
\to \Z \ltimes H_1(M_K,\L)/(t^k-1) \to U(k)$ vanishing on $0 \times P_k$
and lying in
$ P_k^{irr,met}(\pi_1(M_K))$ we get $\eta(M_K,\a)=0$.
\end{theorem}

\begin{proof}
Let $D$ be a slice disk.
 Let $$P_k:=\ker\{H_1(M_K,\L)/(t^k-1)\to H_1(N_D,\L)/(t^k-1)\},$$
this is a metabolizer for the \lp by
proposition \ref{proplinkingformmetabolic}.
 The theorem now follows  from lemma    \ref{lemmaextntond} and
 theorem \ref{letschepkcor}.
\end{proof}

In the following  we will show that some
eta-invariants of slice knots vanish for non-prime power dimensional
irreducible representations.
 Let $\a_1 \in R_{k_1}(G), \a_2 \in R_{k_2}(G)$, then we can form
the tensor product $\a_1 \otimes \a_2 \in R_{k_1k_2}(G)$.

  \begin{proposition} \label{tensorofrep}
If $k_1,k_2$ are coprime, then
\[ \a_{(k_1,z_1,\chi_1)} \otimes \a_{(k_2,z_2,\chi_2)}
\cong  \a_{(k_1k_2,z_1z_2,\chi_1 \chi_2)}. \]
If furthermore $\a_{(k_1,z_1,\chi_1)}$, $\a_{(k_2,z_2,\chi_2)}$
are irreducible, then
$\a_{(k_1,z_1,\chi_1)} \otimes \a_{(k_2,z_2,\chi_2)}$
is irreducible as well.
\end{proposition}

\begin{proof}
 Denote by $e_{11},\dots,e_{k_11}$ and $e_{12},\dots,e_{k_22}$
the canonical bases of
$\C^{k_1}$ and $\C^{k_2}$. Set
$f_i:=e_{i \modd k_1,1} \otimes e_{i \modd k_2,2}$
for $i=0,\dots,k_1k_2-1$.
The $f_i$'s are distinct, therefore $\{f_i\}_{
i=0,\dots,k_1k_2-1}$ form a  basis for $\C^{k_1}\otimes \C^{k_2}$.
One can easily see that $\a_1\otimes \a_2$ with respect to this basis 
is just $\a_{(k_1k_2,z_1z_2,\chi_1\chi_2)}$.

The last statement follows from the observation that if
$\chi_1\chi_2\co H/(t^{k_1k_2}-1)$ factors through $H/(t^k-1)$ for some
$k<k_1k_2$, then one of the $\chi_i\co H\to H/(t^{k_i}-1)$
factors through $H/(t^k-1)$ for some $k<k_i$.
\end{proof}

For a prime number $p$ and a $\L$-module $H$
denote by $P_{k,p}^{irr}(\Z \ltimes H)$
the set of representations $\a_{(z,\chi)}$ in $P_k^{irr}(\Z \ltimes H)$
where $\chi$ factors through a $p$-group.
Define $P_{k,p}^{irr,met}(\pi_1(M_K)):=P_{k,p}^{irr}(\Z \ltimes H_1(M_K,\L))$.

\begin{theorem} \label{mainthm2}
Let $K$ be a slice knot, $k_1,\dots,k_r$ pairwise coprime prime powers, then there
exist  metabolizers $P_{k_i} \subset TH_1(M_{k_i}), i=1,\dots,r$ for the linking pairings, such that
for any prime number $p$ and any choice of
irreducible representations 
$\a_i\co\pi_1(M_K) \to \Z \ltimes H_1(M_K,\L)/(t^{k_i}-1) \to U(k_i)$ vanishing on $0
\times P_{k_i}$ and lying in
$ P_{k_i,p}^{irr,met}(\pi_1(M_K))$ we get
$\eta(M_K,\a_1 \otimes \dots \otimes \a_r)=0$.
\end{theorem}

\begin{proof}
Let $D$ be a slice disk and  let $P_{k_i}:=\ker\{H_1(M_K,\L)/(t^{k_i}-1)\to H_1(N_D,\L)/(t^{k_i}-1)\}$.
All the representations $\a_1,\dots,\a_r$
extend to metabelian representations of $N_D$, hence
$\a_1 \otimes \dots \otimes \a_r$ also extends to a metabelian representation of $N_D$.
Write $\a=\a_{(z_i,\chi_i)}$, then
$\a_1 \otimes \dots \otimes \a_r=
\a_{(z_1\cdot \dots z_r,\chi_1 \cdot \dots \cdot \chi_r)}$
since the $k_i$ are pairwise coprime.
This shows that $\a_1\otimes \dots \otimes \a_r \in P_{k_1\cdot \dots
\cdot k_r,p}^{irr,met}(\pi_1(M_K))$,
therefore $\eta(M_K,\a_1 \otimes
\dots \otimes \a_r)=0$ by theorem  \ref{letschepkcor}.
\end{proof}

\begin{remark}
Theorem  \ref{letschepkcor} holds in fact for locally flat slice disks,
and therefore theorem \ref{mainthm2} holds in fact for topologically slice knots.
Indeed, let $D$
be a topological slice disk for $K$,
i. e. an embedding $D\subset D^4$
such that
$\partial(D)=K$ and  such that $D$ is locally flat, i.e.\   there exists an embedding $f\co D\times D^2\to D^4$ which
extends  the embedding $D\subset D^4$. We write again $N_D:=\overline{D^4\sm f(D\times D^2)}$.
Let $\a \in P_k^{irr,met}(\pi_1(M_K))$ which extends to a representation
$\b \in P_k^{irr,met}(\pi_1(N_D))$.
We can conclude that
$\sign_{\b}(N_D)-k \, \sign(N_D)=0$. Since $N_D$ is in general not a smooth manifold
we can not appeal to  theorem \ref{apsthm} to conclude that $\eta(M_K,\a)=0$.

Now consider the Kirby--Siebenmann invariant $ks(N_D)$ of $N_D$. By \cite[p.\ 10]{R96} we have
$ks(N_D)=\frac{1}{8}\sign(N_D)-\mu(M_K)$ where $\mu(M_K)$ denotes the Rochlin invariant.
Note that $\mu(M_K)$  equals the Arf invariant for $K$, which vanishes for topologically slice knots.
Therefore $ks(N_D)=0$, by \cite[p.\ 125]{FQ90} there exists an $r$ such that $N_D\# rS^2\times S^2$
is smooth. Hence we get
\[\ba{rcl} \eta(M_K,\a)&=&\sign_{\b}(N_D\# rS^2\times S^2)-k \, \sign(N_D\# rS^2\times S^2)\\&=&\sign_{\b}(N_D)-k \, \sign(N_D)=0,\ea\]
where the second to last equality follows from Novikov signature additivity, and the observation that
\[ \sign_{\b}( rS^2\times S^2)=k \, \sign( rS^2\times S^2)\]
since $\b$ is a trivial $k$--dimensional representation.
\end{remark}

Proposition
\ref{tensorofrep} shows that $\a_1\otimes \dots \otimes \a_r$
is irreducible since $\gcd(k_1,\dots,k_r)=1$, i.e.\
the theorem shows that certain non-prime-power dimensional irreducible
eta invariants
vanish for slice knots. Letsche \cite{L00} pointed out the fact that non prime power
dimensional representations can give sliceness (ribbonness) obstructions. In his thesis \cite{F03d}
the author shows that in fact Letsche's non prime
power dimensional representations are tensor products of prime power
dimensional representations.

We say that a knot $K$ has zero slice--eta--obstruction (SE--obstruction)
if the conclusion of theorem \ref{mainthm} holds for all prime powers $k$,
and $K$ has zero slice--tensor--eta--obstruction (STE--obstruction)
if the conclusion of theorem \ref{mainthm2}\break holds for all pairwise coprime prime powers $k_1,\dots,k_r$.

\begin{question}
Are there   examples of knots which
have zero SE--obstruction but non
zero STE--obstruction?
\end{question}

\begin{remark}
It is easy to find examples of a knot $K$ and one-dimensional representations $\a,\b$ such that
$\eta(M_K,\a)=0$ and
$\eta(M_K,\b)=0$ but $\eta(M_K,\a\otimes \b)\ne 0$. This shows that in general
$\eta(M_K,\a\otimes \b)$ is not determined
by $\eta(M_K,\a)$ and $\eta(M_K,\b)$.
\end{remark}

Note that if $TH_1(M_k)=0$, then $R_k^{irr, met}(\pi_1(M_K))=\varnothing$,
therefore
theorem \ref{mainthm2} only gives a non-trivial sliceness obstruction if
$TH_1(M_k)\ne 0$ for some prime power $k$.
This is not
always the case, in fact Livingston  proved the following theorem.

\begin{theorem}{\rm\cite[thm.\ 1.2]{L02}}\label{thmprimepowercovernonzero}\qua
Let $K$ be a knot. There exists a prime power $k$ with  $TH_1(M_k)\ne 0$
if and only if
$\Delta_K(t)$ has a non-trivial irreducible factor that is not an $n$-cyclotomic polynomial with $n$
divisible by three distinct primes.
\end{theorem}

In section  \ref{subsectionexamples},
we will use this theorem to show that there exists a knot $K$ with
$H_1(L_k)=0$ for all prime powers $k$, but
$H_1(L_6)\ne 0$. This shows that $\pi_1(M_K)$ has irreducible
$U(6)$-representations, but no unitary irreducible representations of
prime power dimensions. 
In particular not all representations
are  
tensor products of prime power dimensional representations.\\

\section{Casson-Gordon obstruction}

\subsection{The Casson-Gordon obstruction to a knot being slice}
\label{sectioncg1}
We first recall the definition of the Casson-Gordon obstructions
(cf.\ \cite{CG86}).  
For $m$ a number denote by $C_m \subset S^1$ the unique
cyclic subgroup of order $m$.
For a surjective character $\chi\co H_1(M_k) \to H_1(L_k) \to C_m$,
set
$F_{\chi}:=\Q(e^{2\pi i/m})$.

Since $\Omega_3(\Z \times C_m)=H_3(\Z \times C_m)$
is torsion (cf.\ \cite{CF64}) there exists a 4--manifold $V_k$ and maps $\eps'\co\pi_1(V_k)\to \Z$, $\chi'\co\pi_1(V_k)\to C_m$  such
that 
$\partial(V_k,\eps' \times \chi')=r(M_k,\eps \times \chi)$ for some $r\in \N$. To simplify notation we will denote the maps
$\eps'$ and $\chi'$ on $\pi_1(V_k)$ by $\eps$ and $\chi$ as well.

 The
(surjective) map
$\eps \times \chi\co\pi_1(V_k) \to \Z \times
C_m$ defines a $(\Z \times C_m)$-cover  $\ti{V}_{\infty}$ of $V$.
Then $H_2(C_*(\ti{V}_{\infty}))$
and $F_{\chi}(t)$ have a canonical $\Z[\Z \times C_m]$-module structure
and we can form
$H_2(C_*(\ti{V}_{\infty})
\otimes_{\Z[\Z \times C_m]} F_{\chi}(t))
=\co H_*(V_k,F_{\chi}(t))$.
Since
$F_{\chi}(t)$ is flat over $\Z[\Z \times C_m]$
 by Maschke's theorem (cf.\ \cite{L93})  we get
\[ H_*(V_k,F_{\chi}(t))=H_*(C_*(\ti{V}_{\infty})
 \otimes_{\Z[\Z \times C_m]} F_{\chi}(t))\cong H_*(\ti{V}_{\infty})
\otimes_{\Z[\Z
\times C_m]} F_{\chi}(t),\]
which is a free $F_{\chi}(t)$-module.
If $\chi$ is a character of prime power order, then
the $F_{\chi}(t)$-valued intersection pairing
on $H_2(V_k,F_{\chi}(t))$ is non-singular
(cf.\ \cite[p.\ 190]{CG86}) 
and therefore defines
an element
$t(V_k) \in L_0(F_{\chi}(t))$,
the Witt group of non--singular hermitian forms over $F_{\chi}(t)$ (cf.\ \cite{L93}, \cite{R98}).   
 Denote the image of
the ordinary intersection pairing on $H_2(V_k)$
under the map $L_0(\C) \to L_0(F_{\chi}(t))$ by $t_0(V_k)$.

\CG \cite{CG86} show that 
if $\chi\co H_1(L_k)\to C_m$ is a character with $m$ a prime power, then
\[ \tau(K,\chi):=\frac{1}{r}(t(V_k)-t_0(V_k)) \in L_0(F_{\chi}(t)) \otimes_{\Z} \Q \]
is a well-defined invariant of $(M_k,\eps \times \chi)$,
i.e.\ independent of the choice of $V_k$.
Furthermore they prove the following theorem. 

\begin{theorem}{\rm\cite[p.\ 192]{CG86}}\label{thmcgslice} \qua
Let $k$ be  a prime power. If $K \subset S^3$ is slice
then there exists a metabolizer $P_k$ for the linking pairing,
such that for any $\chi\co TH_1(M_k)\to C_m$, $m$ a
prime power, with
$\chi(P_k) \equiv 0$ we get $\tau(K,\chi)=0$.
\end{theorem}

\subsection{Interpretation of  Casson-Gordon invariants as eta invariants of $M_k$}
\label{sectioncgasetaofmk}
For a number field $F\subset \C$ we consider $F$ with complex involution and $F(t)$ with the
involution given by complex involution and $\bar{t}=t^{-1}$.
For $z\in S^1$ transcendental and $\tau \in L_0(F(t))$ we can consider
$\tau(z)\in L_0(\C)$. Note that $\sign\co L_0(\C)\to \Z$ defines an isomorphism.

Let $K$ be a knot, $k$ any number, $m$ a prime power and
$\chi\co H_1(L_k)\to C_m$ a character and
$(V_k^4,\eps \times \chi\co\pi_1(V_k)\to \Z \times C_m)$ such that
$\partial(V_k,\eps \times \chi)=r(M_k,\eps \times \chi)$.

For a character $\chi\co H_1(L_k) \to C_m$ define characters $\chi^j$
by setting $\chi^j(v):=\chi(v)^j$,
for $z \in S^1$ define
\[ \ba{rcccl} \b_{(z,\chi^j)}\co \pi_1(M_k)& \to&H_1(M_k)=\Z\oplus H_1(L_k)& \to &S^1=U(1) \\
 && (n,v) &\mapsto &z^n \chi^j(v). \ea \]

\begin{proposition} \label{propetamkiscg}
If $z\in S^1$ transcendental, then
\[ \sign(\tau(K,\chi)(z))=\eta(M_K,\b_{(z,\chi^1)}).\]
\end{proposition}

\proof
  Let $z \in S^1$ transcendental, define $\t\co\Z \times C_m \to S^1$ by $\t(n,y):=z^ny$.
Then
\[ \partial(V_k,\t\circ (\eps \times \chi))=r(M_k,\t \circ(\eps \times \chi))=r(M_k,\b_{(z,\chi^1)}). \]
We view $\C$ as a $\Z[\Z\times C_m]$-module
via $\t\circ (\eps \times \chi)$  and
$\C$ as an
$F_{\chi}(t)$ module via evaluating $t$ to $z$.
Note that both modules are flat by Maschke's theorem, hence
\[ \ba{rcl}
H_2^{\b}(V_k,\C)&=&H_2(C_*(\ti{V}_{\infty})) \otimes_{\Z[\Z\times C_m]}\C\\
&=&
(H_2(C_*(\ti{V}_{\infty})) \otimes_{\Z[\Z\times C_m]}F_{\chi}(t))\otimes_{F_{\chi}(t)} \C \\
&=& H_2(V_k,F_{\chi}(t)) \otimes_{F_{\chi}(t)} \C. \ea
\]
This also defines an isometry between the forms, i.e.\ $\sign_{\b}(V_k)=\sign(V_k(z))$.
This shows that
 $$ \ba{rcl} r\eta(M_k,\b(z,\chi^1)) &=&\sign_{\b}(V_k)-\sign(V_k)=\\
&=&\sign(t(V_k)(z))-\sign(t_0(V_k))=r\sign(\tau(K,\chi)(z)).\ea  \eqno{\qed}$$

\begin{remark}
The eta invariant
carries potentially more information than
the function $z \mapsto \sign(\tau(K,\chi)(z))$,
since for non-transcendental $z \in S^1$
the
number $\tau_z(K,\chi)$ is not defined, whereas
$\eta(M_k,\b_{(z,\chi_1)})$ is still defined.
For example the $U(1)$-signatures for slice
knots are zero outside the set of singularities, but
the eta invariant at the singularities contains information
about knots being doubly slice (cf.\ section \ref{sectionu1rep}).
\end{remark}

\begin{proposition} \label{propcgequalsetaformk}
Let $K$ be a knot, $k$ any number, $m$ a prime power and
$\chi\co H_1(L_k)\to C_m$ a character, then the following are equivalent.
\bn
\item $\tau(K,\chi) =0 \in L_0(F_{\chi}(t)) \otimes_{\Z} \Q$.
\item $\sign(\rho(\tau(K,\chi))(z))=0 \in \Q$ for  all transcendental 
$z{\in} S^1, \rho \in
\Gal(F_{\chi},\Q)$.
\item $ \eta(M_k,\b_{(z,\chi^j)})=0 \in \Z \mbox{ for all } (j,m)=1, \mbox{ all transcendental
} z\in S^1.$
\en
\end{proposition}

\proof
The equivalence of (1) and (2) is a purely algebraic statement,
which is shown in a separate paper (cf.\ \cite{F03}) using results of
Ranicki's \cite{R98}.
The equivalence of (2) and (3) follows from proposition \ref{propetamkiscg}
and the observation
that if $\rho \in \Gal(F_{\chi},\Q)$ sends
$e^{2\pi i/m}$ to $e^{2 \pi ij/m}$ for some $(j,m)=1$, then
$\rho(\tau(K,\chi))=\tau(K,\chi^j)$ and hence
$$ \sign(\rho(\tau(K,\chi))(z))=\sign(\tau(K,\chi^j)(z))=
\eta(M_k,\b_{(z,\chi^j)}).\eqno{\qed}$$

\subsection{Interpretation of
Casson-Gordon invariants as eta invariants of $M_K$}
The goal is to prove a version of proposition \ref{propcgequalsetaformk}
with eta invariants of $M_K$
instead of eta invariants of $M_k$.

\begin{proposition}  \label{propmuk}
Let $K$ be a knot, $z\in S^1$ and $\chi\co H_1(M_k)\to H_1(L_k) \to C_m$ a character.
Let $\b:=\b_{(z,\chi)}\co\pi_1(M_k) \to U(1)$ and
$\a=\a_{(z,\chi)}\co\pi_1(M_K) \to U(k)$, then
\[ \eta(M_K,\a)-\eta(M_k,\b)=\sum_{j=1}^k \s_{e^{2\pi ij/k}}(K). \]
If $K$ is algebraically slice and  $H_1(L_k)$ is finite, then 
$\eta(M_K,\a)=\eta(M_k,\b)$.
\end{proposition}

\begin{proof}
In \cite{F03b} we show that if $M_G\to M$ is a $G$--cover and $\a_G\co\pi_1(M_G)\to U(1)$ is
a representation
then
 \[ \eta(M_G,\a_G)=\eta(M,\a)-\eta(M,\a(G)), \]
where $\a(G)\co\pi_1(M)\to G\to U(\C[G])=U(\C^{|G|})$ is given by left multiplication
and $\a\co\pi_1(M)\to U(k)$ is the (induced) representation given by
 \[ \ba{rcl}
\a\co\pi_1(M) &\to& \Aut(\C[ \pi_1(M)]\otimes_{\C[ \pi_1(M_G)]} \C)\\
a &\mapsto &(p\otimes v \mapsto ap \otimes v). \ea
 \]
In our case $G=\Z/k$ and one can easily see that $\a=\a_{(z,\chi)}$.
Since $\Z/k$ is abelian it follows that
$\a(G)=\bigoplus_{i=1}^k \a_i$ where
$\a_i\co\pi_1(M_K)\to U(1)$ is given by $\a_j(z):=e^{2\pi ij/k}$. The proposition now follows from
lemmas \ref{propofcycliccov}  and \ref{propu1sign}.
\end{proof}

We say that a knot $K \subset S^3$ has zero Casson-Gordon obstruction if for any
prime power $k$ there exists a metabolizer $P_k \subset
TH_1(M_k)$ for the linking pairing such that for any
prime power $m$ and
$\chi\co TH_1(M_k)\to C_m$ with
$\chi(P_k) \equiv 0$ we get $\tau(K,\chi)=0 \in L_0(F_{\chi}(t)) \otimes \Q$.

The following is an immediate consequence of propositions  \ref{propcgequalsetaformk} and    \ref{propmuk}.

\begin{theorem} \label{thmetacg}
Let $K$ be an algebraically slice knot. Then $K$ has zero SE--obstruction
if and only if $K$ has
zero Casson-Gordon obstruction.
\end{theorem}

\section{The Cochran-Orr-Teichner-sliceness obstruction} \label{sectioncot}

\subsection{The Cochran-Orr-Teichner-sliceness filtration}

We give a short introduction to the sliceness filtration introduced by Cochran, Orr and Teichner \cite{COT03}.
For a manifold $W$ denote by $W^{(n)}$ the cover corresponding to $\pi_1(W)^{(n)}$. Denote the
equivariant intersection form
\[ H_2(W^{(n)})\times H_2(W^{(n)}) \to \Z [\pi_1(W)/  \pi_1(W)^{(n)}] \]  
by $\lambda_n$, and the
self-intersection form by $\mu_n$. 
An $(n)$-Lagrangian is a submodule $L \subset H_2(W^{(n)})$ on which $\lambda_n$ and $\mu_n$ vanish
and which maps onto a Lagrangian of $\lambda_0 \co H_2(W)\times H_2(W)\to \Z$.

\begin{definition}\cite[def.\ 8.5]{COT03}  
A knot $K$ is called $(n)$-solvable if $M_K$ bounds a spin 4-manifold $W$ such that $H_1(M_K)\to H_1(W)$ is an
isomorphism and such that $W$ admits two dual $(n)$-Lagrangians. This means that $\l_n$ pairs the
two Lagrangians non-singularly and that the projections freely generate $H_2(W)$.

A knot $K$ is called $(n.5)$-solvable if $M_K$ bounds a spin 4-manifold $W$ such that $H_1(M_K)\to H_1(W)$ is an
isomorphism and such that $W$ admits an $(n)$-Lagrangian and a dual $(n+1)$-Lagrangian.
 
$W$ is called an $(n)$-solution respectively an $(n.5)$-solution for $K$.
 \end{definition}

\begin{remark}$\phantom{9}$
\bn
\item The size of an $(n)$-Lagrangian depends only on the size of $H_2(W)$,
in particular if $K$ is slice, $D$ a slice disk, then $N_D$
is an
$(n)$-solution for
$K$ for all
$n$, since
 $H_2(N_D)=0$.
\item By the naturality of covering spaces and homology with twisted coefficients it follows that
if $K$ is $(h)$-solvable, then $K$ is $(k)$-solvable for all $k<h$.
\en
\end{remark}

\begin{theorem}$\phantom{9}$
\bn
\item The following are equivalent.
\bn
\item[\rm(a)] $K$  is (0)--solvable.
\item[\rm(b)] $\Arf(K)=0$.
\item[\rm(b)] $\Delta_K(t)\equiv \pm 1 \mod 8$.
\en
\item A knot $K$ is (0.5)--solvable if and only if $K$ is algebraically slice.
\item If $K$ is (1.5)--solvable then $K$ is algebraically slice and $K$ has vanishing Casson-Gordon invariants.
\item There exist algebraically slice knots which have
zero Casson-Gordon invariants but are not $(1.5)$--solvable.
\item There exist $(2.0)$--solvable knots which are not slice.
\en
\end{theorem}

Statements (1), (2), 3) and (5) are shown in \cite{COT03}.
Taehee Kim \cite{K02} showed that there exist $(1.0)$--solvable knots which have zero Casson-Gordon invariants, but
 are not $(1.5)$--solvable (cf.\ also proposition \ref{example4}).

For any $n\in \N$ Tim Cochran and Peter Teicher have examples (unpublished)
 of knots that are
 $(n)$--solvable but not $(n.5)$--solvable, in particular are not slice.
 It is unknown whether for $n\in \N_0$ every $(n.5)$--solvable knot is $(n+1)$--solvable.

\subsection{$L^2$--eta invariants as sliceness-obstructions}
In this section we very  quickly summarize some   $L^2$--eta invariant theory.

Let $M^3$ be a smooth manifold and $\varphi\co\pi_1(M)\to G$ a homomorphism,
then Cheeger and Gromov \cite{CG85} defined an invariant
$\eta^{(2)}(M,\varphi)\in \R$, the (reduced)
$L^2$--eta invariant. When it is clear which homomorphism we mean, we will write
$\eta^{(2)}(M,G)$ for $\eta^{(2)}(M,\varphi)$.

\begin{remark}
If  
$\partial(W,\psi)=(M^3,\varphi)$,
then (cf.\ \cite[lemma 5.9 and remark 5.10]{COT03}, \cite{LS03})
\[ \eta^{(2)}(M,\varphi)=\sign^{(2)}(W,\psi)-\sign(W), \]
where $\sign^{(2)}(W,\psi)$  denotes Atiyah's
\cite{A76}   $L^2$--signature.
\end{remark}

Let $\Q\L:=\Q[t,t^{-1}]$.

\begin{theorem}{\rm {\cite{COT03}}}\qua
 \bn
\item If $K$ is $(0.5)$--solvable, then $\lteta(M_K,\Z)=0$.
\item If $K$ is $(1.5)$--solvable, then there exists a metabolizer $P_{\Q}\subset H_1(M_K,\Q\L)$ for the Blanchfield
pairing
\[ \l_{Bl,\Q}\co H_1(M_K,\Q\L) \times H_1(M_K,\Q\L) \to \Q(t)/\Q[t,t^{-1}] \]
such that for all $x\in P_{\Q}$ we get $\lteta(M_K,\b_x)=0$, where $\b_x$ denotes the map
\[\ba{rl} \pi_1(M_K)\to \Z \ltimes H_1(M_K,\L) &\to \Z \ltimes H_1(M_K,\Q\L)\\
&\xrightarrow{id \times \l_{Bl,\Q}(x,-)}\Z \ltimes \Q(t)/\Q[t,t^{-1}].\ea\]
\en
\end{theorem}

We say that a knot $K$ has zero abelian $L^2$--eta invariant sliceness obstruction  if
$\lteta(M_K,\Z)=0$ and $K$  has zero metabelian $L^2$--eta invariant  sliceness obstruction
if there exists a metabolizer $P_{\Q}\subset H_1(M_K,\Q\L)$ for $\l_{Bl,\Q}$
such that for all $x\in P_{\Q}$ we get $\lteta(M_K,\b_x)=0$.
Note that if $K$ has zero metabelian $L^2$--eta invariant   then it is easy to see that the integral Blanchfield form is metabolic
as well and hence  $K$ is algebraically slice by proposition \ref{keartonthm}.

\section{Examples} \label{sectionex}
In this section we will construct 
\bn
\item a knot which has zero abelian $L^2$--eta invariant, but is not algebraically slice,
\item a $(1.0)$--solvable knot which has zero metabelian $L^2$--eta invariant, but non-zero SE--obstruction,
\item a knot which has zero $STE$--obstruction but non-zero metabelian $L^2$--eta invariant
(following Taehee Kim \cite{K02}).
\en
The idea in examples (2) and (3) as well as in the examples of section \ref{examplerib} is to start out with a slice knot $K$
and make `slight' changes via
a satellite construction.  The change in the eta invariants
can be computed explicitly.

\subsection{Satellite knots}
 Let $K,C \in S^3$ be knots. Let $A\subset S^3\sm K$ a simple closed curve, unknotted in $S^3$, note that $S^3\sm N(A) $ is a torus.
Let
$\varphi\co\partial(N(A)) \to \partial(N(C))$ be a diffeomorphism which sends a meridian of $A$ to a
longitude of $C$ and a longitude of $A$ to a meridian of $C$. The space
\[ \overline{S^3 \sm N(A)} \cup_{\varphi} \overline{S^3\sm N(C)} \]
is a 3--sphere and the image of $K$ is denoted by $S=S(K,C,A)$. We say $S$ is the satellite knot
with companion $C$, orbit $K$ and axis $A$.
Note that this construction is equivalent to replacing a tubular neighborhood
of $C$ by the torus knot $K\subset \overline{S^3\sm N(A)}$.

\begin{proposition} \label{satelliteknotslice}
If $C$ is slice, then for any $K$ and $A$ the satellite knot $S(K,C,A)$  is concordant to $K$.
If $K\subset S^1\times D^2\subset S^3$ and $C$ are  ribbon, then any $S(K,C,A)$ is in fact ribbon.

\end{proposition}

\begin{proof}
Let $K \subset   S^1\times D^2 \subset S^3$ and let $C$ be a slice knot.
Let $\phi\co S^1 \times I \to S^3 \times I$ be a null-concordance for $C$, i.e.\ $\phi(S^1 \times 0)=C$ and $\phi(S^1
\times 1)$
is the unknot. We can extend this to a map $\phi\co S^1\times D^2 \times I \to S^3 \times I$ such that $\phi\co S^1 \times
D^2\times
0$ is the zero-framing for $C$.

Pick a diffeomorphism $f\co\overline{S^3\sm N(A)}\to S^1\times D^2$ such that the meridian and longitude of $A$
get sent to the longitude and meridian of $S^1\times 0$.
Now consider
\[ \psi\co S^1 \times I \to K\times I \hookrightarrow \overline{S^3\sm N(A)}\times I\xrightarrow{f\times id}  S^1\times
D^2\times I\xrightarrow{\phi} S^3
\times I.\] Note that $\phi\co S^1 \times D^2\times 1$ is a zero framing for the unknot, since linking numbers are
concordance invariants and $\phi\co S^1 \times D^2\times 0$ is the zero framing for $C$. This shows that $\psi\co S^1\times 1
\to S^3
\times 1$ gives the  satellite knot of the unknot with orbit
$K$, i.e.\ $K$ itself. Therefore $\psi$ gives a
concordance between $S=\psi(S^1 \times 0)$ and
$K=\psi(S^1 \times 1)$.

Now assume that $K,C$ are ribbon. Then we can find a concordance $\phi$ which has no minima under
the projection $S^1\times [0,1] \to S^3\times [0,1]\to [0,1]$. It is clear that $\psi$ also has no minima,
capping off with a ribbon disk for $K$ we get a disk  bounding $S$ with no minima, i.e.\ $S$ is ribbon.
\end{proof}

\begin{proposition}{\rm\cite[p.\ 8]{COT04}}\label{thmsolvabilityofsat}\qua 
Let $K$ be an $(n.0)$--solvable knot,
$C$ a $(0)$--solvable knot, $A \subset S^3\sm K$ such that $A$ is the unknot in $S^3$
and
$[A]\in
\pi_1(S^3\sm K)^{(n)}$. Then $S=S(K,C,A)$ is $(n)$--solvable.
\end{proposition}

\subsection{Eta invariants   of  satellite knots}
Let $S$ be a satellite knot with companion $C$, orbit $K$ and axis $A$. 
\begin{proposition}{\rm\cite[p.\ 337]{L84}}\label{propsksamea}\qua
If $A\in \pi_1(S^3\sm K)^{(1)}$ then the inclusions
\[ S^3\sm N(K)\hookrightarrow S^3 \sm N(K)\sm N(C) \hookleftarrow
(\overline{S^3 \sm N(A)} \cup_{\varphi} \overline{S^3\sm N(C)})\sm N(K) = S^3 \sm N(S) \]
induce isometries of the Blanchfield pairing and the linking pairing of $S$ and $K$ for any cover such that
$H_1(L_{K,k})$ is finite.
\end{proposition}

This lemma shows in particular that we can
identify the set of characters on $H_1(L_{K,k})$
with the set of characters on $H_1(L_{S,k})$.
The following corollary follows immediately from \cite{T73}.
\begin{corollary}\label{corsksamea}
If $S=S(K,C,A)$ as above, then $S$ and $C$ have $S$--equivalent Seifert matrices.
\end{corollary}

Let $S=S(K,C,A)$ be a satellite knot with $A\in \pi_1(S^3\sm K)^{(1)}$.
Let $k$ be any number such that $H_1(L_{K,k})$ is finite.
The curve $A \subset S^3\sm N(K)$ is null-homologous
 and therefore lifts to simple closed curves $\ti{A}_1,\dots,\ti{A}_k \in L_{K,k}$.
\begin{theorem} \label{etaforsatellite}
 \label{coretaforsatellite}
Let $z \in S^1$ and $\chi\co H_1(L_{K,k})=H_1(L_{S,k}) \to S^1$
a  character. Then
\[ \eta(M_S,\a^S_{(\chi,z)})=\eta(M_K,\a^K_{(\chi,z)}) + \sum_{i=1}^{k}
\eta(M_C,\a_i), \]
where $\a_i$ denotes the representation $\pi_1(M_C) \to U(1)$
given by $g \mapsto \chi(\ti{A}_i)^{\eps(g)}$.
\end{theorem}

\begin{proof}
 Litherland \cite{L84} proved a general statement
how to compute the Casson-Gordon invariant of $S$ in terms of the
Casson-Gordon invariant of $K$ and the basic invariants  of $C$. Translating the proof into the language of
eta invariants shows that
\[ \eta(M_{S,k},\b^S_{(\chi,z)})=\eta(M_{K,k},\b^K_{(\chi,z)}) + \sum_{i=1}^{k}
\eta(M_C,\a_i). \]
For full details we refer to \cite{F03d}.
Casson-Gordon invariants only make sense
when  $k$ and the order of $\chi$ are prime powers,
but the proof can be used to show the above statement about eta invariants.

 By corollary \ref{corsksamea} the knots $K$ and $S$
have $S$--equivalent Seifert matrices, in particular the twisted
signatures are the same
(cf.\ proposition \ref{propu1sign}). The theorem now follows immediately from
proposition \ref{propmuk}.
\end{proof}

Iterating the satellite construction we can generalize the theorem as follows (cf.\ \cite[p.\ 405]{L02}).
Let $K\in S^3$ be a knot and $A_1,\dots,A_s \in S^3\sm K$ be simple closed curves which form
the unlink in $S^3$ and such that $A_i=0 \in \pi_1(S^3\sm K)^{(1)}$. Let $C_1,\dots,C_s$ be knots. Then we can
inductively form satellite knots by setting
$S_0:=K$ and $S_i$ the satellite formed with orbit $S_{i-1}$, companion $C_i$ and axis $A_i$.
Note that $A_i \in \pi_1(S^3\sm S_{i-1})^{(1)}$. We write
\[ S_i=\co S(K,C_1,\dots,C_i,A_1,\dots,A_i). \]

\begin{theorem} \label{etaforiteratedsatellite}
Let $S:=S(K,C_1,\dots,C_s,A_1,\dots,A_s)$ as above.
Let $k$ be any number  such that $H_1(L_{K,k})$ is finite, $z \in S^1$ and $\chi\co H_1(L_{S,k}) \to S^1$ a  character,
denote the corresponding character $H_1(L_{K,k}) \to S^1$ by $\chi$ as well.
Then
\[ \eta(M_{S},\a^S_{(\chi,z)})=\eta(M_{K},\a^K_{(\chi,z)}) + \sum_{j=1}^s \sum_{i=1}^{k}
\eta(M_{C_j},\a_{ij}). \]
Here $\a_{ij}$ denotes the representation $\pi_1(M_{C_j}) \to U(1)$
given by $g {\mapsto} \chi((\ti{A}_j)_{i})^{\eps(g)}$, where
$(\ti{A}_j)_1,\dots,(\ti{A}_j)_k$ denote the lifts of $A_j$ to $L_{S,k}$.
\end{theorem}

\subsection{Computation of  $L^2$--eta invariants}
The following proposition makes it possible to compute $L^2$--eta invariants in many cases.

\begin{proposition}[\cite{COT03}, \cite{F03c}]
Let $K$ be a knot, then
\[
  \eta^{(2)}(M_K,\Z)=\int_{S^1} \s_{z}(K).  \]
\end{proposition}

We  quote a theorem by \COT on the computation of $L^2$--eta invariants for
satellite knots. Recall that for $x\in H_1(M_,\L)$ we defined
\[\ba{rl} \b_x\co \pi_1(M_K)\to \Z \ltimes H_1(M_K,\L) &
\to \Z \ltimes H_1(M_K,\Q\L)\\
&\xrightarrow{id \times \l_{Bl,\Q}(x,-)}\Z \ltimes \Q(t)/\Q[t,t^{-1}].\ea\]

\begin{theorem}[{\cite[p.\ 8]{COT04} \cite[prop.\ 5.3]{K02}}]\label{thml2forsat}
Let $S=S(K,C,A)$ with $A\in \pi_1(S^3\sm K)^{(1)}$.
Let $x\in H_1(M_S,\L)=H_1(M_K,\L)$, then
\[ \eta^{(2)}(M_S,\b_x)=
\left\{ \ba{ll}
\eta^{(2)}(M_K,\b_x)+\eta^{(2)}(M_C,\Z) & \mbox{ if }\b_x(A)\ne 0\\
 \eta^{(2)}(M_K,\b_x)  & \mbox{ if }\b_x(A)= 0  \ea \right.
\]
where  $A$ is considered as an element in $H_1(M_K,\L)$.
\end{theorem}

In \cite{F03c} we give a different approach to the computation of metabelian $L^2$--eta invariants
based on an approximation theorem by L\"uck and Schick \cite{LS01}.

\subsection{Examples} \label{subsectionexamples}

$$ B_1:= \bp 0 &0&1&1 \\ 0&0&0&1\\  1&1&0&1\\  0&1&0&0 \ep\leqno{\rm Let}$$
this Seifert matrix is obviously metabolic. The Alexander polynomial is
$\Delta_{B_1}(t)$\break $=(t^2-t+1)^2$.
The signature function $z \mapsto \s_z(B_1)$ is zero outside of the set of \zeros of the Alexander polynomial since the
form is metabolic. The \zeros are $e^{2\pi i/6}, e^{2\pi 5i/6}$ and at both points the signature is
$-1$.
Let
\[ B_2:= \bp 1&1 \\ 0&1 \ep. \]
Then $\Delta_{B_2}(t)=t^2-t+1$, and 
\[ \s_{e^{2\pi i t}}(B_2)= \left\{ 
\ba{ll} 2 & \mbox{ for }  t\in (\frac{1}{6},\frac{5}{6})  \\
0 & \mbox{ for }  t\in [0,\frac{1}{6}) \cup  
  (\frac{5}{6},1] \ea
    \right. \]
 Finally let 
\[ B_3:= \bp 1&0&0&0&1&0\\ -1&1&0&1&0&1\\ 0&0&0&1&1&1 \\0&1&0&1&0&1\\ 1&0&1&0&1&1\\ 0&1&1&1&0&1\ep \]
Then $\Delta_{B_3}(t)=\Phi_{14}(t)=1-t+t^2-t^3+t^4-t^5+t^6$, and 
\[ \s_{e^{2\pi i t}}(B_3)= \left\{
\ba{ll} 2 & \mbox{ for }  t\in (\frac{1}{14},\frac{3}{14}) \cup  
(\frac{5}{14},\frac{9}{14}) \cup (\frac{11}{14},\frac{13}{14}) \\
0 & \mbox{ for }  t\in [0,\frac{1}{14}) \cup  
(\frac{3}{14},\frac{5}{14}) \cup (\frac{9}{14},\frac{11}{14})\cup (\frac{13}{14},1] \ea
    \right. \]

\begin{proposition}[Example 1] \label{propexzerol2}
There exists a $(0)$--solvable knot $K$ with zero abelian $L^2$--eta invariant   but which
is not algebraically torsion.
\end{proposition}

\begin{proof}
Recall that for a knot $\Arf(K)=0$ if and only if $\Delta_K(-1) \equiv \pm 1 \mod 8$.
By \cite{L69} we can find a knot $K$  with Seifert matrix $7B_3\oplus -6B_2$,
then $\arf(K)=7\arf(B_3)-6\arf(B_2)=0$, using the above calculations we get
 $\s_z(K)=2$ for $z=e^{2\pi i k/5},k=1,2,3,4$ and
$\int_{S^1}\s_z(K)=0$.
This shows that $K$ has all the required properties.
\end{proof}

\begin{proposition}[Example 2]  \label{propex2}
There exists a   $(1.0)$--solvable knot $K$
with zero metabelian $L^2$--eta invariant, but non-zero SE-invariants.
\end{proposition}

\begin{proof}
Let $p(t)=-2t+5-2t^{-1}$. 
By Kearton \cite{K73} there exists a knot $K_5\cong S^5\subset S^7$ such that
its Blanchfield pairing is isomorphic to
\[ \ba{rcl} \L/p(t)^2\times \L/p(t)^2&\to &S^{-1}\L/\L \\
 (a,b) &\mapsto & \bar{a}p(t)^{-2}b \ea \]
Let $A$ be a Seifert matrix for $K_5$, and $K\subset S^3$ be a slice knot with Seifert matrix $A$.
Since the Blanchfield pairing is determined by $A$, the Blanchfield pairing of $K$ is isomorphic to the Blanchfield pairing of $K_5$.
A computation using lemma  \ref{propofcycliccov}
shows that $|H_1(L_{K,4})|=225$.

Let $N$ be a non--negative integer greater or equal than
$X:=\mbox{max}\{\eta(M_K,\a) | \a \in R_4(\pi_1(M_K)) \}$.
Note that $R_4(\pi_1(M_K))$ is compact since $\pi_1(M_K)$ is finitely generated.
Using of results of Levine \cite[p.\ 92]{L94} one can show that $X$ is finite
(cf.\ \cite{F03c} for details).

Let $F$ be a Seifert surface for $K$. We can view $F$ as a disk with $2g$ 1-handles attached. 
Then the meridians $A_1,\dots,A_{2g} \in S^3 \sm F$ of the handles 
 form the unlink in $S^3$ and   the corresponding homology classes give a basis for $H_1(S^3
\sm F)$. Note that $A_i \in\pi_1(S^3\sm K)^{(1)}$.
Denote the knot of the proof of proposition \ref{propexzerol2} by $D$
and let $C=(N+1)\cdot D$, and form the
iterated satellite knot
\[ S:=S(K,C,\dots,C,A_1,\dots,A_{2g}). \]
We claim that the satellite knot $S$ satisfies the conditions stated in the proposition.
 $S$ is $(1.0)$--solvable by theorem \ref{thmsolvabilityofsat} and has zero 
metabelian $L^2$--eta invariant 
by theorem \ref{thml2forsat} since $K$ is slice, $\int_{S^1}\s_z(C)=0$ by construction of $C$
and since $K$, and therefore also $S$, has a unique metabolizer for the Blanchfield pairing.

We have to show that for all $P_4 \subset H_1(L_{S,4})$ with
$P_4=P_4^{\perp}$ with respect to the \lp $\l_{S,4}$,
we can find a non-zero character $\chi\co H_1(L_{S,4})\to S^1$ of prime power order,
vanishing on $P_4$, such that for one transcendental $z$
we get $\eta(M_S,\a_{(z,\chi)})\ne 0$.

Let $P_4$ be a metabolizer and $\chi\co H_1(L_{S,4})\to S^1$
a non-trivial character of order 5,
vanishing on $P_4$. Denote the corresponding character on $H_1(L_{K,4})$ by $\chi$ as well.
For any $z\in S^1$ we get by  corollary \ref{coretaforsatellite}
\[ \eta(M_{S},\a^S_{(\chi,z)})=\eta(M_{K},\a^K_{(\chi,z)}) + \sum_{j=1}^{2g} \sum_{i=1}^{4}
\eta(M_{C},\a_{ij}), \]
where $\a_{ij}$ denotes the representation $\pi_1(M_{C}) \to U(1)$
given by $g \mapsto \chi((\ti{A}_j)_{i})^{\eps(g)}$ and $(\ti{A}_j)_i$ denotes the $i^{th}$ lift of $A_j$ to $L_{K,k}$.
By definition of $N$ and by proposition \ref{propu1sign} we get
\[ \ba{rcl}\eta(M_{S},\a^S_{(\chi,z)})& \geq& -N+\sum_{j=1}^{2g} \sum_{i=1}^{4}
\eta(M_{C},\a_{ij})\\
&=&-N+\sum_{j=1}^{2g} \sum_{i=1}^{4}\s_{\chi((\ti{A}_j)_i)}(C)\\
&=&-N+\sum_{j=1}^{2g} \sum_{i=1}^{4}(N+1)\s_{\chi((\ti{A}_j)_i)}(B_2)\ea \]
Note that $\eta(M_{C},\a_{ij})\geq 0$ for all $i,j$ since
$\s_{e^{2\pi i j/5}}(C)\geq 0$  for $j=0,\dots,4$.   The lifts $(\ti{A}_{j})_i$ are easily seen to generate
$H_1(L_{K,4})$, hence $
\chi((\ti{A}_j)_{i}) \ne 1$ for at least one $(i,j)$ since $\chi$ is
non-trivial.
But $\s_w(C)=2(N+1)$ for
$w=e^{2\pi ij/5}, j=1,2,3,4$. It follows that
$\eta(M_{S},\a^S_{(\chi,z)}) \geq -N+2(N+1)> 0$ for all $z$.
\end{proof}

For completeness sake we add the following example which was discovered by 
Taehee Kim \cite{K02}. 
\begin{proposition}[Example 3] \label{example4}
There exists a knot $S$ which is algebraically slice, $(1.0)$--solvable, has zero STE--obstruction but
non-zero metabelian $L^2$--eta invariant.
\end{proposition}

\begin{proof}
Denote by 
$\Phi_{30}(t)=1+t-t^3-t^4-t^5+t^7+t^8$ the minimal polynomial of $e^{2\pi i/30}$.
As in the proof of proposition \ref{propex2} there exists a ribbon knot $K$ such that
the Blanchfield pairing is isomorphic to
\[ \ba{rcl} \L/\Phi_{30}(t)^2\times \L/\Phi_{30}(t)^2&\to &S^{-1}\L/\L \\
 (a,b) &\mapsto & \bar{a}\Phi_{30}(t)^{-2}b \ea \]
An explicit example of such a knot is given by Taehee Kim \cite[Section 2]{K02}. 
Note that $K$ has a unique metabolizer $P$ for the Blanchfield pairing.  
Furthermore
$H_1(L_{K,k})=0$ for all prime powers $k$
by theorem \ref{thmprimepowercovernonzero}.
Let $C$ be a knot with Seifert matrix $B_1$ and $A\in \pi_1(S^3\sm K)^{(1)}$ unknotted in $S^3$
and such that $\b_x(A)\ne 0$ for some $x\in P$.

By corollary   \ref{corsksamea}    the knot
$S$ is algebraically slice, since $K$ is algebraically slice. $S$ is $(1.0)$--solvable
by proposition \ref{thmsolvabilityofsat} since $\arf(C)=0$.

Since $H_1(L_{S,k}){=}H_1(L_{K,k}) {=}0$ for all prime powers $k$,
we get $R_k^{irr,met}(\pi_1(M_K))$ $=\emptyset$ for all prime powers $k$,
hence $S$ has zero STE--obstruction.
By theorem \ref{thml2forsat} $S$ has non-zero metabelian $L^2$--eta invariant   since  
$\b_x(S)=\b_x(K)+\int_{S^1}\s_z(C)=\frac{4}{3}$ and $P$ is the unique
metabolizer.
\end{proof}

\section{Ribbon knots, doubly slice knots and the obstructions of Gilmer and Letsche}

\subsection{Obstructions to a knot being ribbon} \label{sectionrib}

\begin{figure}[ht!]
\begin{center}
\includegraphics[scale=0.25]{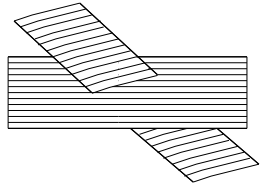}
\caption{Immersed ribbon disk.}
\label{ribbonsing}
\end{center}
\end{figure}

An immersed ribbon disk is an immersion $D\to S^3$ which bounds $K$
such that the singularities are only of the type as in figure
\ref{ribbonsing}, i.e.\ the self--intersection lies completely in the interior
of one of the two sheets involved. The singularities can be resolved in $D^4$
to give an embedded slice disk. Slice disks which are isotopic to such
disks are called ribbon disks.
If a knot has a ribbon disk we say that the knot is ribbon. Note that this
is a purely 3--dimensional definition.

It is a longstanding conjecture of Fox (cf.\ \cite[problem 25]{F61}) that all
slice knots are ribbon. In theorem \ref{mainribbonthm}
we give a condition for a knot to be ribbon which is ostensibly stronger
than the corresponding condition (theorem \ref{mainthm2})
for a knot to be slice. It is an intriguing question whether one can use these
results to disprove Fox's conjecture.

\begin{proposition}[{\cite[lemma 3.1]{G81} \cite[lemma 2.1]{K75b}}]\label{pi1mnsurjective}
If $K$ is ribbon and $D \subset D^4$ a ribbon disk,  then the maps
\[ \ba{rcl} i_*\co\pi_1(S^3 \sm K) &\to & \pi_1(D^4\sm D)\\
 H_1(M_K,\L) &\to &H_1(N_D,\L)
\ea \]
are surjective.
\end{proposition}

\begin{proposition}  \label{propribbonmetab}
Assume that $K$ is ribbon, $D$ a ribbon disk, then $H_1(N_D,\L)$ is $\Z$-torsion free, in particular
$P:=\ker\{
H_1(M_K,\L)
\to H_1(N_D,\L)\}$ is a metabolizer for $\l_{Bl}$.
\end{proposition}

\begin{proof}
According to \cite[thm.\ 2.1 and prop.\ 2.4]{L77} there exists a short exact
sequence
\[ 0 \to \Ext_{\L}^2(H_1(N,M,\L)) \to \overline{H_1(N,\L)} \to
\Ext_{\L}^1(H_2(N,M,\L)) \to 0.\]
Here $ \overline{H_1(N,\L)}$ denotes $H_1(N,\L)$ with involuted
$\L$-module
structure, i.e.\ $t\cdot v:=t^{-1}v$.
Furthermore $\Ext_{\L}^1(H_2(N,M,\L))$ is $\Z$-torsion
free (cf.\
\cite[prop.\ 3.2]{L77}).
In order to show that
 $H_1(N,\L)$ is $\Z$-torsion free
it is therefore
enough to show that $H_1(N,M,\L)=0$.
Consider the exact sequence
\[ H_1(M,\L) \to H_1(N,\L) \to H_1(N,M,\L) \to H_0(M,\L) \to H_0(N,\L) \to
0.\]
The last map is an
isomorphism.
By proposition \ref{pi1mnsurjective} the first map is surjective. It follows
that
$H_1(N,M,\L)=0$.

The second part follows immediately from theorem \ref{firstpropmetabl}.\end{proof}

\begin{theorem} \label{mainribbonthm}
Let $K \subset S^3$ be a ribbon knot.
Then there exists a metabolizer $P \subset
H_1(M_K,\L)$ such that
 for any $\a \in P_k^{irr}(\pi_1(M_K))$ vanishing on $0 \times P$
we get $\eta(M_k,\a)=0$.
\end{theorem}

\begin{proof}
Let $P:=\ker\{ H_1(M_K,\L) \to H_1(N_D,\L)\}$ where $D$ is a ribbon
disk for $K$. Then $P=P^{\perp}$ by
proposition \ref{propribbonmetab}.
Let $\a \in P_k(\pi_1(M))$ which vanishes on $0 \times P$, then $\a$ extends to a metabelian
representation of $\pi_1(N_D)$, hence $\eta(M_K,\a_{(z,\chi)})=0$
by lemma    \ref{lemmaextntond} and theorem \ref{letschepkcor}.
\end{proof}

We say that $K$ has zero eta invariant ribbonness obstruction if the conclusion holds for $K$.

\begin{remark}
One can show that if $P$ is a metabolizer for $\l_{Bl}$ and
$k$ is such that $H_1(L_k)$ is finite, then $P_k:=\pi_k(P) \subset TH_1(M_k)=H_1(M,\L)/(t^k-1)$
is a metabolizer for $\l_{L,k}$. In particular if
$K$ is and $k$ any number such that $H_1(L_k)$ is finite then
there exists a metabolizer $P_k$ for the linking pairing  such that for all
$\chi\co TH_1(M_k)
\to S^1$ of prime power order, vanishing on
$P_k$, and for all transcendental $z\in S^1$ we get $\eta(M_K,\a_{(z,\chi)})=0$.
 
Comparing  this result with theorem \ref{mainthm2}  we see that the ribbon obstruction  
is stronger in two respects. When $K$ is ribbon
\bn
\item we can find metabolizers $P_k$ which all lift to the same metabolizer of the \blpp,
\item the representations for non-prime power dimensions don't have to be tensor products.
\en
We will make use of this in proposition \ref{example3}.
\end{remark}

\begin{remark}
\bn
\item
Note that the only fact we used was that for a ribbon disk $H_1(M_K,\L)\to H_1(N_D,\L)$ is surjective.
\item
\CG \cite[p.\ 154]{CG86} prove a ribbon--obstruction theorem which does not require
the character to be of prime power order, but has a strong restriction on the
fundamental group of $\pi_1(N_D)$.
\item
In \cite{F03c} we show that if a knot $K$ has zero  eta invariant ribbonness obstruction, then $K$ has
in particular zero metabelian $L^2$--eta invariant sliceness obstruction.
In proposition \ref{propex3} we will see that the converse is not true.
\en
\end{remark}

\subsection{Obstructions to a knot being \ds} \label{sectiondoublyslice}
A knot $K\subset S^3$ is called \ds (or doubly null-concordant) if there exists an unknotted smooth
two-sphere $S
\subset S^4$ such that
$S\cap S^3=K$.
It follows from the Schoenflies theorem that a \ds knot is in particular slice. Fox \cite{F61} posed the question
which slice knots are doubly slice. Doubly slice knots have been studied by
Sumners \cite{S71}, Levine \cite{L89},
Ruberman \cite{R83} and Taehee Kim \cite{K02}.

We say that knot $K$ is algebraically doubly slice if $K$ has a Seifert matrix of the form
$\bp 0 & B \\ C&0 \ep$ where $B,C$ are square matrices of the same size.
 Sumners \cite{S71} showed
that if $K\subset S^3$ is doubly slice, then $K$ is algebraically doubly slice.
This result can be used to show that many slice knots are not doubly slice.

We prove the following  new doubly sliceness obstruction theorem.

\begin{theorem} \label{maindoublyslicethm}
Let $K \subset S^3$ be a \ds knot.
Then there exist metabolizers $P_1, P_2 \subset
H_1(M_K,\L)$ for the Blanchfield pairing such that
\bn
\item $H_1(M_K,\L)=P_1 \oplus P_2$,
\item
 for any $\a \in P_k^{irr}(\Z \ltimes H_1(M_K,\L))$ vanishing on $0 \times P_i$, i=1,2,
we get $\eta(M_K,\a)=0$.
\en
\end{theorem}

\begin{proof}
Let $S\subset S^4$ be an unknotted two-sphere such that $S \cap S^3=K$. Intersecting
$S$ with
$\{ (x_1,\dots,x_5) \in \R^5 | x_5\geq 0\}$ and
$\{ (x_1,\dots,x_5) \in \R^5 | x_5\leq 0\}$ we can write
write $S=D^2_1\cup_K D^2_2$ and $S^4=D_1^4\cup_{S^3} D_2^4$.
Let $N_i=D^4_i\sm D^2_i$, then $N_1\cap N_2=S^3\sm K$ and $N_1\cup N_2=S^4\sm S$.
From  the Mayer-Vietoris sequence
we get
\[ H_1(M_K,\L)=H_1(N_1,\L)\oplus H_1(N_2,\L)\]
since $H_1(M_K,\L)=H_1(S^3\sm K,\L)$
and $H_1(S^4\sm S,\L)=0$ since $S$ is trivial.

Now let $P_i:=\ker\{ H_1(M,\L) \to H_1(N_i,\L)\}$,
the proof concludes as the proof of the main ribbon obstruction theorem \ref{mainribbonthm}
since $H_1(M_K,\L)\to H_1(N_i,\L)$ is surjective for $i=1,2$.
\end{proof}

\begin{remark}
\bn
\item The proof of theorem \ref{maindoublyslicethm} shows in particular  that
if $K$ is \ds we can find a slice disk $D$ such that $H_1(N_D,\L)$ is $\Z$-torsion free.
\item Comparing theorem \ref{maindoublyslicethm} with theorems \ref{mainthm2}
and \ref{mainribbonthm} we see that \ds knots have
zero (doubly) ribbon obstruction.
\en
\end{remark}

\begin{question}
Using the notation of the proof we get from the van Kampen theorem that
for a \ds knot
\[ \Z =\pi_1(D_1^4\sm D^2_1)*_{\pi_1(S^3\sm K)} \pi_1(D_2^4\sm D^2_2).\]
Can we conclude that $\pi_1(S^3\sm K)\to \pi_1(D_i^4\sm D^2_i)$ is surjective for at least one $i$?
If yes, then this would show that a \ds knot is in fact homotopically ribbon.
One can go further and ask whether doubly slice knots are in fact ribbon or doubly ribbon.
The first part of the question is of course a weaker version of the famous `slice equals ribbon' conjecture.
\end{question}

\begin{remark}
Taehee Kim \cite{K03}, \cite{K04} introduced the notion of $(n,m)$-solva\-bility ($n,m\in \frac{1}{2}\N$).
He used
$L^2$-eta invariants to find highly non--trivial examples of non doubly slice knots,
in particular for each $n\in \N$ he found examples of knots which are $(n,n)$--solvable
but not $(n+1,n+1)$--solvable.
\end{remark}

\subsection{The Gilmer--obstruction} \label{sectiongilmer}
In sections \ref{sectiongilmer} and \ref{sectionletsche}
we quickly recall the obstructions of Gilmer
\cite{G83}, \cite{G93} and Letsche \cite{L00}.
The papers claim to define sliceness obstructions.
While investigating the precise connection to our obstructions we found
that unfortunately the proofs in their paper contain gaps.
We show that Gilmer's and Letsche's results give ribbon obstructions
and we give a precise description of the problems that must be solved in
order to show that they are really sliceness obstructions.

Let $K$ be a knot, $F$ a Seifert surface. Pick a basis $a_1,\dots,a_{2g}$
for $H_1(F)$, denote by $A$ the corresponding Seifert matrix.
Let $\G:=(A^t-A)^{-1}A^t$ and $k$ such that $H_1(L_k)$ is finite.
Define $\varphi_k\co H_1(F) \to H_1(F)$ to be
the endomorphism given by $\G^k-(\G-1)^k$ and
define $B^k \subset H_1(F,\Q/\Z)$ to be the
kernel of $\varphi_k \otimes \Q/\Z$. For a prime number $p$
define $B^k_p$ to be the $p$-primary part of $B^k$.
Let $Y$ be $S^3$ cut along $F$ and denote by $\a_1,\dots,\a_{2g} \in H_1(Y)$ the dual basis with respect to Alexander
duality, i.e.\ $\lk(a_i,\a_j)=\delta_{ij}$.

Pick a lift of  $Y$ to $L_k$.
Denote the lifts of $\a_i$ by $\ti{\a}_i$. These generate $H_1(L_k)$, in fact
\[ H_1(L_k) = (\Oplus \Z\ti{\a}_i) /\G_k^t. \]
We therefore get a well--defined isomorphism
\[ \ba{rcl} B^k &\to & H^1(L_k,\Q/\Z)=\hom(H_1(L_k),\Q/\Z) \\
    \sum r_ia_i &\mapsto & (\ti{\a}_j \mapsto r_j)
    \mbox{ where }r_j \in \Q/\Z. \ea \]
Different lifts of $Y$ to $L_k$ give different characters, but  the associated 
Casson-Gordon invariants agree. Hence we get a well-defined Casson-Gordon invariant $\tau(K,\chi) \in
L_0(F_{\chi}(t))\otimes \Q$ for $\chi \in B^k$.

We say that a knot $K$ has zero Gilmer obstruction for a Seifert surface $F$ and a prime $p$ if 
there exists a
metabolizer $H$ for the Seifert pairing
on $H_1(F)$ such that for all  prime 
powers $k$, $\tau(K,B^k_p \cap (H \otimes \Q/\Z))=0$.
More precisely, for a character
$\chi \in B^k_p \cap (H \otimes \Q/\Z)$
we get $\tau(K,\chi)=0 \in
L_0(F_{\chi}(t))$.

Gilmer's theorem (cf.\ \cite[p.\ 5]{G93}) claims that a slice knot has
zero Gilmer obstruction for all Seifert surfaces and all primes. Unfortunately the proof has a gap.
On page 6, the statement that $H \otimes \Q/\Z$ is the kernel of $\mu_*$
(Gilmer's notation in the paper) is not necessarily true
since tensoring with $\Q/\Z$ is not exact.
This becomes a problem if $H_1(R)$ is not torsion free,
where $R$ denotes a 3--manifold which bounds the union of $F$ with a slice disk $D$.
One can show that the question whether $H_1(R)$ is torsion free is closely related
to the question
whether $\ker\{H_1(M_K,\L)\to H_1(N_D,\L)\}$ is a metabolizer for the Blanchfield
pairing.

Furthermore the proof of the cancellation lemma 5  has a gap as well,
namely on the second to last line.
Note that the same problem appears in Gilmer's earlier paper \cite{G83}.\\

Going carefully through the proof of Gilmer's theorem, or the equivalent version in terms of eta invariants in \cite{F03d} one
can see the following.

\begin{theorem}
Let $K$ be a slice knot and $F$ a Seifert surface. Then $K$ has zero Gilmer obstruction for $F$ and all but finitely many  primes $p$.
\end{theorem}

In \cite{F03d} we show that if a knot satisfies the vanishing condition on eta invariants
of theorem \ref{mainribbonthm} then the Gilmer obstruction vanishes for all Seifert
surfaces
and all primes. In particular this shows
the following theorem.

\begin{theorem} \label{gilmerstheorem}
Let $K$ be a ribbon knot or a doubly slice knot, then $K$ has zero Gilmer obstruction for all $F$ and all primes.
\end{theorem}

\subsection{The  Letsche--obstruction} \label{sectionletsche}

For $x \in H_1(M_K,\L)$ we define the map
\[\ba{rl}  \a_x\co\pi_1(M_K) \to \pi_1(M_K)/\pi_1(M_K)^{(2)}& \cong \Z \ltimes
H_1(M_K,\L)\\& \xrightarrow{id \times \l_{Bl}(x,-)} \Z \ltimes S^{-1}\L/\L. \ea  \]
We say that a knot $K$ has zero  Letsche obstruction
if
there exists a metabolizer $P\subset H_1(M_K,\L)$ such that
 for any $k$, any $x\in P$ and any $\t  \in R_k(\Z \ltimes S^{-1}\L/\L)$ such that
$\t\circ\a_x \in P_k^{irr}(\pi_1(M_K))$ we get $\eta(M_K,\t \circ \a_x)=0$.
Letsche (cf.\ \cite[p.\ 313]{L00}) claims that every slice knot has zero Letsche obstruction.
Unfortunately the statement of the last paragraph of the proof
of lemma 2.21 is incorrect since maps $\chi_1,\chi_2$ to
$S^1$ whose $n^{th}$ powers $\chi_1^n,\chi_2^n$ agree for some $n>1$, don't have to be identical,
i.e.\ $\chi_1\ne \chi_2$.
In particular given an abelian group $P$ and a map
$\chi\co nP\to S^1$ for some $n>1$ there is no canonical way to extend
$\chi$ to $P$.

The situation in Letsche's paper is as follows. Let $D$ be a slice disk for $K$,
then by theorem \ref{keartonthm} $P:=\ker\{H_1(M_K,\L)\to FH_1(N_D,\L)\}$
is a metabolizer, but only characters which vanish
on $Q:=\ker\{H_1(M_K,\L)\to H_1(N_D,\L)\}$ extend over the slice disk complement.
Note that $Q=nP$ for some $n\in \N$.
If $x\in P$ and $n\ne 1$, then $\a_x$ does not vanish on the metabolizer
$P$ but only on $Q=nP$.
Letsche's attempt to get around this problem introduced the above mentioned problem.

On the other hand if $Q=P$, i.e.\
if $\ker\{H_1(M_K,\L) \to H_1(N_D,\L)\}$ is a metabolizer for the
Blanchfield pairing, then the Letsche obstruction vanishes.
In particular we get the following weaker statement.

\begin{theorem}  \label{thmletsche}
Let $K \subset S^3$ be a slice knot, $D$ a slice disk.
If  $\ker\{H_1(M_K,\L)$ $\to H_1(N_D,\L)\}$ is a metabolizer for the
Blanchfield pairing, then $K$ has zero Letsche obstruction.
In particular ribbon knots and doubly--slice knots have zero Letsche obstruction.
\end{theorem}

\begin{remark}
Note that there is no restriction on the dimensions of the representations. In \cite{F03d} we show that all unitary irreducible representations
of $\Z \ltimes S^{-1}\L/\L$ are in fact tensor products of unitary representations of prime power dimensions.
\end{remark}

We give a complete proof, which differs somewhat from Letsche's original proof.

\begin{proof}
Let $x \in P$.
Considering the long exact sequence we see that
$x=\partial(w)$ for some $w \in H_2(N_D,M_K,\L)$.
Letsche \cite{L00} showed that
in fact $H_2(N_D,M_K,\L)=\tor_{\L}H_2(N_D,M_K,\L)$ and that there exists
a \blp
\[ \l_{Bl,N_D}\co\tor_{\L}H_2(N_D,M_K,\L) \times \tor_{\L}H_1(N_D,\L) \to S^{-1}\L/\L \]
such that $\l_{Bl}(x,y)=\l_{Bl,N_D}(w,i_*(y))$ for $y\in H_1(M_K,\L)$.
We get a commutative diagram (cf.\ \cite[cor.\ 2.9]{L00})
\[ \ba{ccccccc}
\pi_1(M_K) &\to & \Z {\ltimes} H_1(M_K,\L) &\xrightarrow{id \times \l_{Bl}(x,-)} & \Z {\ltimes} S^{-1}\L/\L
&\xrightarrow{\t} & U(k) \\
\downarrow &&\downarrow &&\parallel &&\parallel \\
\pi_1(N_D) &\to & \Z {\ltimes} \tor_{\L}H_1(N_D,\L) &\xrightarrow{id \times \l_{Bl,N_D}(w,-)} & \Z {\ltimes} S^{-1}\L/\L
&\xrightarrow{\t}&
 U(k). \ea \]
This shows that $\t \circ \a_x$ extends over $\pi_1(N_D)$.
The first part of the theorem now follows from theorem \ref{letschepkcor}.
The second part follows from theorem \ref{firstpropmetabl}, proposition  \ref{propribbonmetab}
and the remark after theorem \ref{maindoublyslicethm}.
\end{proof}

\subsection{More Examples} \label{examplerib}
In \cite{F03c} we show that if for a knot the metabelian eta invariant ribbon obstruction vanishes,
then all the abelian and metabelian  sliceness obstructions vanish.
The following example shows that furthermore the ribbon obstruction theorem is in fact stronger than
these obstructions, considered as ribbonness obstructions.

 \begin{proposition}[Example 4]\label{example3} \label{propex3}
There exists a knot $S$ which is algebraically slice, $(1.0)$--solvable, which has zero STE--obstruction and zero
metabelian $L^2$--eta invariant but which does not satisfy the
condition for  theorem \ref{mainribbonthm}, i.e.\
$S$ is not ribbon.
\end{proposition}

\begin{proof}
Let $K$ be the ribbon knot of proposition \ref{example4}. Recall that $\Delta_K(t)=\Phi_{30}(t)^2$ and that $K$ has a unique metabolizer $P$.
In particular
$H_1(L_{K,k})=0$ for all prime powers $k$ and a computation using
lemma \ref{propofcycliccov} shows that
$|H_1(L_{K,6})|=625$.

Now let $C$ be the knot of the proof of proposition \ref{propexzerol2},
recall that $\int_{S^1} \s_z(C)=0$.
Let $P_6:=\pi_6(P)\subset H_1(M_{K,6})$ be the projection of $P$.
Let $A$ be a simple closed curve in $S^3\sm K$, unknotted in $S^3$,
such that $A\in \pi_1(S^3\sm K)^{(1)}$,
which lifts to a simple closed curve $\ti{A}$
in the
$6$-fold cover which presents a non--trivial element of order 5 in $H_1(M_{K,6})/P_6$.

We claim that the satellite knot $S:=S(K,C,A)$ satisfies the conditions stated in the proposition.
The proof of proposition  \ref{example4} shows that $S=S(K,C,A)$ is algebraically slice, is $(1.0)$--solvable, has zero STE--obstruction
and zero metabelian $L^2$--eta invariant obstruction.

As remarked above, the Blanchfield pairing of $S$ has a unique metabolizer $P$.
Let
\[ \chi\co H_1(M_K,\L)\to H_1(M_K,\L)/(t^6-1)\to H_1(L_{S,6})\to S^1\] be
a non-trivial character of order 5,
vanishing on $P_6\subset H_1(L_{S,6})$ such that $\chi(\ti{A})\ne 1\in S^1$.
By corollary \ref{coretaforsatellite}
\[ \eta(M_S,\a^S_{(\chi,z)})=\eta(M_K,\a^K_{(\chi,z)}) + \sum_{i=1}^{k}
\eta(M_C,\a_i), \]
where $\a_i$ denotes the representation $\pi_1(M_C) \to U(1)$
given by $g \mapsto \chi(\ti{A}_i)^{\eps(g)}$.
The first term is zero since $K$ is ribbon and $P$ is the unique metabolizer
of the Blanchfield pairing (cf.\ theorem \ref{mainribbonthm}). 
The second  term is non-zero since $\chi(\ti{A}_i)\ne 1$ for at least one $i$
and by the properties of $C$.
This  shows that 
$\eta(M_{S},\a^S_{(\chi,z)})\ne 0$, i.e.\ $S$ is not ribbon by theorem \ref{mainribbonthm}.
\end{proof}

\begin{remark}
\bn
\item
The knot $S$ constructed in the proof has in fact the extra property that
the Gilmer and the Letsche obstruction vanish. Indeed, in unpublished work
we show that the Gilmer
(respectively Letsche) obstruction
is equivalent to the vanishing of certain eta invariants corresponding to irreducible
prime power dimensional representations (respectively tensor products thereof). By construction
the knot $S$ has no irreducible prime power dimensional representations.
\item
It is in fact possible to construct an example of a knot with all the above properties such
that no multiple of $S$ is ribbon.
\en
\end{remark}

Note that it is not known whether the example given in the above proposition is slice or not.

\begin{proposition}[Example 5]\label{example5}
There exists a ribbon knot $S$ with the following property. There exists a prime power $k$ such that
for all metabolizers $P$ we can find $\a \in P_5^{irr,met}(\pi_1(M_K))$ with $\a(0\times P)=0$ but such that
 $\eta(M_K,\a) \ne 0$.
\end{proposition}

Note that we do not restrict ourselves to representations $\a_{(z,\chi)}$, with $\chi$ of prime power order.
The proposition shows that theorems \ref{mainthm} and \ref{mainribbonthm}
can't be strengthened to include all
non prime power characters. The example shows that the set $P_k^{irr,met}(\pi_1(M_K))$ is in a sense maximal,
i.e.\ that the prime power condition on the characters is indeed necessary.

\begin{proof}
Consider
$\Delta(t)=f(t)f(t^{-1})$ where $f(t)=4-3t+2t^2+4t^3-7t^4+t^5+2t^6-3t^7+t^8$.
 Terasaka \cite{T59} shows that any polynomial of
the form $f(t)f(t^{-1})$ can be realized by a slice knot. But this shows that $\Delta(t)$ can be realized by a
metabolic Seifert matrix, hence there exists a ribbon knot $K$
with $\Delta_K(t)=\Delta(t)$.
A computation  shows that $H_1(L_{K,5})=1296=36^2$.
Let $N$ be an integer greater than $\max\{\eta(M_K,\a)|\a \in R_5(\pi_1(M_K))\}$.
 As in the proof of proposition  \ref{propex2} we see that such an $N$ exists.

Let $\ti{A}_1,\dots,\ti{A}_s \in H_1(L_{K,5})$ be all elements.
One can find simple
closed curves $A_1,\dots,A_s \subset S^3\sm K$
 such that $A_i\in \pi_1(S^3\sm K)^{(1)}$ and  such that for all $i=1,\dots,s$ the homology class $\ti{A}_i$ is
represented by one of the $k$ lifts of $A_i$ to $L_k$.
Possibly after crossing changes of the representatives one can assume that
$A_1,\dots,A_s $ form in fact the unlink in $S^3$.
Let $C$  be a ribbon knot with Seifert matrix $\oplus_{i=1}^{N+1}B_1$, and form the iterated satellite knot
\[ S:=S(K,C,\dots,C,A_1,\dots,A_s). \]
Note that $S$ is ribbon by proposition \ref{satelliteknotslice}.

Let $P$ be a metabolizer for the Blanchfield pairing and let
\[ P_5:=\pi_5(P)\subset  H_1(M_K,\L)/(t^5-1)=H_1(L_{K,5}).\]
Let $\chi\co H_1(L_5)\to
S^1$ be a non-trivial character of order 6, vanishing on
$P_5$. Then for all transcendental $z\in S^1$ we get
\[ \ba{rcl} \eta(M_{S,5},\a^S_{(\chi,z)})&=&\eta(M_{K,5},\a^K_{(\chi,z)}) + \sum_{j=1}^s \sum_{i=1}^{5}
\eta(M_{C_j},\a_{ij}) \\
&\leq &N + \sum_{j=1}^s \sum_{i=1}^{5}
(N+1) \s_{\chi((\ti{A_j})_i)}(B_1). \ea \]
Since $\chi$ is of order 6 and by construction of $A_1,\dots,A_s$, we can find $(i,j)$ such  that
$\chi((\ti{A_j})_i)=e^{2\pi i/6}$, but recall that $\s_z(B_1)=0$ for all $z$ except for $z=e^{2\pi i/6}, e^{2\pi
5i/6}$ where $\s_z(B_1)=-1$. This shows that $\eta(M_{S},\a^S_{(\chi,z)}) \leq N+(N+1)(-1)=-1$.
 \end{proof}

\begin{remark}
The degree of the Alexander polynomial of $K$ is quite large, but it
was the polynomial of lowest degree we could find with
$|H_1(L_k)|$ being divisible by $6$ for some prime power $k$.
\end{remark}

\Addresses\recd

\end{document}